\documentclass[11pt]{article}
\usepackage{numbysec,epstopdf}
\usepackage{graphics,graphicx,fancyhdr,color}

\hsize \textwidth
\advance \hsize by -\marginparwidth
\evensidemargin \oddsidemargin
\usepackage{amssymb}
\advance\hoffset by 5mm

\newtheorem{thm}[theorem]{Theorem}
\newtheorem{prop}[theorem]{Proposition}
\newtheorem{lemma}[theorem]{Lemma}

\newcommand{\Rm}{{\mathbb R}}
\newcommand{\pdr}[2]{\frac{\partial{#1}}{\partial{#2}}}

\newcommand{\eps}{\varepsilon}

\newcommand{\commentout}[1]{}
\newcommand{\disp}{\displaystyle}
\numberbysection

\begin{document}
\title{The explosion problem in a flow}
\date{}
\author{Henri Berestycki\thanks{EHESS, CAMS,
54 Boulevard Raspail, F - 75006 Paris, France; hb@ehess.fr}\and
Alexander Kiselev\thanks{Department of
Mathematics, University of Wisconsin, Madison, WI 53706; e-mail:
kiselev@math.wisc.edu} \and
Alexei Novikov\thanks{Department of
Mathematics, Pennsylvania State University, University Park, PA
16802; anovikov@math.psu.edu}
\and
Lenya Ryzhik\thanks{Department of Mathematics, University of
Chicago, Chicago, IL 60637, USA; ryzhik@math.uchicago.edu}}
\maketitle
\begin{abstract}
We consider the explosion problem in an incompressible flow introduced
in \cite{BKJS}. We use a novel $L^p-L^\infty$ estimate for elliptic
advection-diffusion problems to show that the explosion threshold
obeys a positive lower bound which is uniform in the advecting
flow. We also identify the flows for which the explosion threshold
tends to infinity as their amplitude grows and obtain an effective
description of the explosion threshold in the strong flow asymptotics
in a two-dimensional one-cell flow.
\end{abstract}

\section{Introduction}

\subsection{The explosion problem}

The explosion problem concerns existence and regularity of
positive solutions of nonlinear elliptic equations of the form
\begin{equation}\label{intro-ellipt}
-\Delta\phi=\lambda g(\phi),
\end{equation}
in a domain $\Omega\subset\Rm^n$ with the Dirichlet boundary
conditions:
$\phi=0$ on the boundary $\partial\Omega$. The
nonlinearity $g(\phi)$ is convex and increasing with $g(0)>0$ and
\begin{equation}\label{gint-infty}
\int_0^\infty\frac{ds}{g(s)}<+\infty.
\end{equation}
Two typical examples to keep in mind are $g(s)=e^s$ and $g(s)=(1+s)^m$
with $m>1$. The positive parameter $\lambda>0$ measures the
non-dimensional strength of the nonlinearity. It has been shown in the
pioneering works of Keener and H. Keller \cite{KK}, Joseph and
Lundgren \cite{JL}, and Crandall and Rabinowitz \cite{CR} that there
exists a critical threshold $\lambda^*>0$ so that (\ref{intro-ellipt})
admits positive solutions for $0<\lambda<\lambda^*$, while no positive
solutions exist for $\lambda>\lambda^*$. The regularity of solutions
at $\lambda=\lambda^*$ is a delicate issue: the linearized problem
was studied by Brezis and Vazquez in \cite{BV} in great detail.
In particular, when the domain is a ball, and
for the exponential and power nonlinearities mentioned above, the
solutions at the
critical value $\lambda^*$ are uniformly bounded in dimensions less or
equal to $N=9$ and $N=10$, respectively, while in higher dimensions
they are unbounded. For more general nonlinearities $g(s)$
and domains $\Omega$, regularity of solutions at $\lambda=\lambda^*$ in
dimensions $N=2,3$ has been established by Nedev \cite{Nedev}, and
more recently in dimension $N=4$ by Cabr\'e \cite{Cabre-06}.

In the present paper we consider the non-selfadjoint elliptic problem
\begin{eqnarray}\label{elliptic}
&&-\Delta\phi+u\cdot\nabla\phi=\lambda g(\phi)
\hbox{ in $\Omega,$}\\
&&\phi=0\hbox{ on $\partial\Omega,$}\nonumber
\end{eqnarray}
with a prescribed incompressible flow $u(x)$ so that $\nabla\cdot
u=0$.  Extension of the aforementioned results to the case when a flow
is present is a natural question in the context of the original
motivation for the study of (\ref{intro-ellipt}) as the explosion
problem \cite{FK,Semenov,ZBLM}. Of particular interest is to understand
how the presence of an underlying flow and its features affect the explosion
limit.

Existence of a critical explosion threshold $\lambda^*(u)$ can be
established as a straightforward generalization of existing
methods. We are mostly interested here in the qualitative
dependence of $\lambda^*(u)$ on the flow $u$ -- whether a flow may
raise or lower the explosion threshold, and in the asymptotic
behavior of $\lambda^*(u)$ in the limit of a strong flow.
Intuitively, a flow improves mixing and interaction with the
boundary -- hence one may expect that an incompressible flow would
always raise the explosion threshold. Somewhat surprisingly, this
was shown not necessarily to be the case in \cite{BKJS}.  More
precisely, Berestycki, Kagan, Joulin and Sivashinsky have
considered in \cite{BKJS}, the problem (\ref{elliptic}) for a
two-dimensional cellular flow and observed numerically that while
the explosion threshold increases for flows oscillating on a small
scale, it may actually decrease if the flow has large scale
variations. This is because such flows may promote creation of hot
spots where the explosion would happen faster than without any
flow. The authors of \cite{BKJS} have also presented a formal
asymptotic analysis and found an effective problem in the limit of
the large flow amplitude. The fully nonlinear problem when the
flow itself satisfies a Navier-Stokes type equation coupled to the
explosion problem for temperature has been studied in
\cite{Volpert-Chaos,JMS,MS} using numerics and formal asymptotics.
Recently, some rigorous results for the behavior of the solutions
to the coupled system in the regime of a strong gravity have been
obtained in \cite{CNR}. Here we derive several qualitative
properties of the explosion threshold $\lambda^*(u)$ in terms of
the geometry and the amplitude of the flow $u$.

\subsection{The main results}

In the following we always assume that $\Omega$ is a smooth bounded domain in $\Rm^n$,
$u(x)$ is a $C^1(\bar\Omega)$ divergence-free flow ($\nabla\cdot u=0$ in $\Omega$).
Our first proposition establishes the direct analog of the classical
results for (\ref{intro-ellipt}) and allows us to define the
critical parameter $\lambda^*(u)$.
\begin{prop}\label{thm1}
There exists $\lambda^*(u)\in(0,\infty)$ such that (i) for every
$0<\lambda<\lambda^*(u)$ the problem (\ref{elliptic}) has a unique
positive classical solution $\phi_\lambda(x)$
such that the principal eigenvalue
$\kappa_1$ of the linearized operator
$M\psi=-\Delta\psi+u\cdot\nabla\psi- \lambda g'(\phi_\lambda)\psi$
is positive; (ii) if (\ref{elliptic}) admits another non-negative
solution $v(x)$ then $v(x)\ge\phi_\lambda$; (iii) the function
$\phi_\lambda(x)$ is increasing in $\lambda$; (iv) there exists no
classical solution of (\ref{elliptic}) for $\lambda>\lambda^*(u)$.
\end{prop}
The proof of this result is very close to that in \cite{CR} -- we
present it below both for the convenience of the reader and since we
will use some of the intermediate steps in what follows. Another
reason to discuss the proofs for $u\not\equiv 0$ is that some of the
basic results in the self-adjoint case $u=0$ rely on the variational
characterization of the principal eigenvalue of the linearized
operator $M$ which we do not have when the flow is present.

The next theorem shows that the possible creation of hot spots
cannot drop the explosion threshold arbitrarily close to zero, no
matter what the incompressible flow $u(x)$ is.
\begin{thm}\label{thm-lowerbd}
For any domain $\Omega$ and nonlinearity $g(\phi)$ there exists
$\lambda_0>0$ so that the critical threshold $\lambda^*(u)$ for
(\ref{elliptic}) satisfies $\lambda^*\ge\lambda_0>0$ for all
incompressible flows $u(x)$ in $\Omega$. The constant $\lambda_0$
depends on $\Omega$ and the function $g$.
\end{thm}
This result does not hold without the restriction that the flow $u(x)$
is incompressible -- we describe in Section \ref{sec:uniform} examples
of flows for which $\lambda^*$ may be as small as one wishes. The proof
of Theorem~\ref{thm-lowerbd} involves the following  uniform $L^p-L^\infty$ bound
for solutions of the Dirichlet problem for elliptic
diffusion-advection problems with the constant independent of the
incompressible flow.
\begin{lemma}\label{lem-lp-linfty}
Let the flow $u(x)$ be divergence-free and let $q(x)$ be the solution of the
elliptic problem
\begin{eqnarray}\label{lbd-lem}
&&-\Delta q+u\cdot\nabla q=f(x)\hbox{ in $\Omega$,} \\
&&q=0\hbox{ on $\partial\Omega$,} \nonumber
\end{eqnarray}
with $f(x)\in L^p(\Omega)$, $p>n/2$. There exists a constant $C(\Omega,n,p)>0$
which depends on $p$ and the domain $\Omega$ but not on the flow $u(x)$,
so that $\|q\|_{L^\infty(\Omega)}\le C\|f\|_{L^p(\Omega)}$.
\end{lemma}
These results raise several interesting open questions. First, can
one identify the optimal constant $C$ in (\ref{lbd-lem})? More
importantly, we would like to pose as an open problem to know
whether in all domains with smooth boundaries a flow realizing the
best constant in Lemma~\ref{lem-lp-linfty} exists. The same
question pertains to the smallest possible explosion threshold in
Theorem~\ref{thm-lowerbd}: does the flow minimizing $\lambda^*(u)$
over all incompressible flows exist in any domain $\Omega$? If so,
what are its geometric characteristics? Numerical simulations
in~\cite{BKJS} indicate that a natural guess that $u=0$ turns out
not to be correct for all domains as an incompressible flow may
create additional hot spots.  More precisely, it has been
numerically computed  in~\cite{BKJS} that in a very long rectangle
the explosion threshold corresponding to a cellular flow with a
certain finite positive amplitude is smaller than that
corresponding to $u=0$. However, it is not clear whether in the
situation when $u=0$ is not a minimizer of $\lambda^*(u)$,   such
a minimizer  exists at all, or if minimizing flows do not exist.
When it exists, how is it determined?

Let us now fix a flow profile $u(x)$ and consider the explosion
problem (\ref{elliptic}) with a strong flow $Au(x)$, with a large flow
amplitude $A\gg 1$:
\begin{eqnarray}\label{elliptic-strong}
&&-\Delta\phi+Au\cdot\nabla\phi=\lambda g(\phi)
\hbox{ in $\Omega$}\\
&&\phi=0\hbox{ on $\partial\Omega$.}\nonumber
\end{eqnarray}
We are interested in the behavior of the explosion threshold $\lambda^*(A)$
for (\ref{elliptic-strong}) in the limit $A\to +\infty$. Let us recall that
a function $\psi\in H^1(\Omega)$ is a first integral of $u$ if $u\cdot\nabla\psi=0$ a.e.
in $\Omega$.
\begin{thm}\label{thm5}
We have $\lambda^*(A)\to+\infty$ as $A\to+\infty$ if and only if $u$
has no non-zero first integrals in $H^1_0(\Omega)$.
\end{thm}
This theorem provides a sharp characterization of the flows capable of
preventing an explosion for an arbitrary $\lambda>0$ provided that the
flow amplitude $A$ is sufficiently large. The proof uses the ideas
from \cite{BHN} and \cite{CKRZ} together with some techniques
of \cite{BCMR}. Not surprisingly, the explosion threshold tends
to infinity as $A\to+\infty$ under the same assumptions as the principal
Dirichlet eigenvalue of the operator $-\Delta +Au\cdot\nabla$ (see \cite{BHN}),
as both quantities measure the effectiveness of the enhancement of the
boundary cooling due to the flow.

Finally, in Section~\ref{sec:cellflow}, we consider the effective
problem for (\ref{elliptic-strong}) in the limit $A\to +\infty$
for the class of two-dimensional cellular flows. We show that in
this limit, the various cells of the flow ``do not talk to each
other''. The main result of that section is
Theorem~\ref{cf-fr-thm}. In particular, when $A\to+\infty$,  the
explosion threshold $\lambda^*(A)$ is close to the explosion
threshold  on the ``largest'' cell in $\Omega$. Moreover, the
explosion threshold for each of the individual flow cells in the
limit  $A\to+\infty$ has an asymptotic description in terms of the
Freidlin problem. We recall that the fast flow asymptotics for the
parabolic reaction-diffusion equations in flows without cells has
been treated by M. Freidlin in \cite{Freidlin}, and our results
for the elliptic explosion problem in the special case of one cell
flows are what one would expect formally from \cite{Freidlin}. The
most interesting and delicate new ingredient is  the independent
behavior of the solutions in various cells when $A\to+\infty$. We
mention that unlike in~\cite{Freidlin}, our proofs are not
probabilistic in nature. Actually, as a by-product of the present
paper, one can use our arguments to recover some of the results
of~\cite{Freidlin} by analytic techniques.

Another natural variable coefficients extension of the classical
results for (\ref{intro-ellipt}) is to allow the nonlinearity to be
spatially dependent -- work in this direction
has been recently done in \cite{GG1,GG2}.

{\bf Acknowledgment.} This work has been carried as the first
author was visiting the University of Chicago. It has been
supported by ASC Flash Center at the University of Chicago. AK was
supported by NSF grant DMS-0314129, AN  was supported by NSF grant
DMS-0604600, LR was supported by NSF grant DMS-0604687.

\section{Existence and basic properties of $\lambda^*$}

\subsection{The  critical parameter}\label{sec:basic}

We begin with the proof of Proposition \ref{thm1} which is well
known for $u=0$. We do not assume in this section that the flow
$u(x)$ is incompressible. As we have mentioned, the proof is very
close to that in \cite{CR}, with some minor modifications.
Let $\mu_1[u]$ and $\eta(x)$ be the
principal eigenvalue and the normalized positive eigenfunction of
the adjoint problem
\begin{eqnarray}\label{elliptic-0mu}
&&-\Delta\eta-\nabla\cdot(u\eta)=\mu_1[u]\eta
\hbox{ in $\Omega$}\\
&&\eta=0\hbox{ on $\partial\Omega$.}\nonumber
\end{eqnarray}
Note that $\mu_1[u]>0$ as the operator  $-\Delta+u\cdot\nabla$
has no zero order term.
\begin{lemma}\label{prop-large-lambda}
The problem (\ref{elliptic}) admits no non-negative classical
solutions for $\lambda>
\mu_1[u]/g'(0)$.
\end{lemma}
{\bf Proof.} Since the function $g(s)$ is convex and $g(0)>0$ we have
$g(s)\ge g'(0)s$. Therefore, any classical solution $\phi_\lambda\ge 0$ of
(\ref{elliptic}) satisfies
\begin{eqnarray}\label{elliptic-012}
&&-\Delta\phi_\lambda+u\cdot\nabla\phi_\lambda\ge \lambda g'(0)\phi_\lambda
\hbox{ in $\Omega,$}\\
&&\phi_\lambda=0\hbox{ on $\partial\Omega$.}\nonumber
\end{eqnarray}
Multiplying (\ref{elliptic-012}) by $\eta$ and integrating by parts
we conclude that
\[
\lambda g'(0)\int \eta\phi_\lambda \le
\int\phi_\lambda[-\Delta\eta-u\cdot\nabla\eta]= \mu_1[u]\int
\eta\phi_\lambda.
\]
It follows that for a positive solution of (\ref{elliptic}) to
exist we must have $\mu_1[u]\ge g'(0)\lambda$ and thus no
non-negative solution of (\ref{elliptic}) exists if
$\lambda>\mu_1[u]/g'(0)$. $\Box$

Next, we show that for a sufficiently small $\lambda>0$ a positive
solution of
 (\ref{elliptic}) exists. Let
$\tau(x)$ be the expected value of the exit-time, solution of
\begin{eqnarray}\label{eq-tau0}
&&-\Delta \tau+u\cdot\nabla \tau=1\hbox{ in $\Omega,$}\\
&&\tau=0\hbox{ on $\partial\Omega,$}\nonumber
\end{eqnarray}
and let
\begin{equation}\label{def-thetau}
\theta_u=\max_{x\in\bar\Omega}\tau(x).
\end{equation}
\begin{lemma}\label{prop-small-lambda}
There exists a constant $C>0$ which depends only on the
nonlinearity $g(s)$ but not on the flow $u(x)$ so that problem
(\ref{elliptic}) admits a minimal non-negative solution
$\phi_\lambda$ for all $\lambda\le C/\theta_u$.
\end{lemma}
The proof is by constructing a super-solution and using it to show
that a positive solution of (\ref{elliptic}) exists. Let us recall
the following fact.
\begin{lemma}\label{lem-super}
Assume that there exists  a smooth function $\bar\phi(x)\ge 0$  satisfying
\begin{eqnarray}\label{eq-barphi}
&&-\Delta\bar\phi+u\cdot\nabla\bar\phi\ge
\lambda g(\bar\phi)\hbox{ in $\Omega$,}\\
&&\bar\phi\ge 0\hbox{ on $\partial\Omega$.}\nonumber
\end{eqnarray}
Then there exists a classical solution $\phi_\lambda$ of (\ref{elliptic})
which is minimal in the sense that for any other non-negative
solution $\psi$
of (\ref{elliptic}) we have $\phi_\lambda(x)\le\psi(x)$ for all
$x\in\Omega$.
\end{lemma}
{\bf Proof.} We construct an approximating sequence $\phi_n(x)$ by
setting $\phi_0(x)=0$ and letting $\phi_{n+1}$ be the smooth
solution of
\begin{eqnarray}\label{elliptic=0n}
&&-\Delta\phi_{n+1}+u\cdot\nabla\phi_{n+1}=\lambda g(\phi_n)
\hbox{ in $\Omega$,}\\
&&\phi_{n+1}=0\hbox{ on $\partial\Omega$.}\nonumber
\end{eqnarray}
The difference $w_1:=\phi_1-\phi_0(=\phi_1)$ satisfies
\begin{eqnarray}\label{elliptic=0nw}
&&-\Delta w_1+u\cdot\nabla w_1=\lambda g(0)\ge 0
\hbox{ in $\Omega$,}\\
&&w_1=0\hbox{ on $\partial\Omega$.}\nonumber
\end{eqnarray}
It follows from the maximum principle that $w_1\ge 0$ and thus
$\phi_1\ge\phi_0$. Similarly, we have for the higher differences
$w_n=\phi_n-\phi_{n-1}$:
\begin{eqnarray}\label{elliptic=0nwn}
&&-\Delta w_n+u\cdot\nabla w_n=\lambda
[g(\phi_{n-1})-g(\phi_{n-2})]\ge 0
\hbox{ in $\Omega$,}\\
&&w_1=0\hbox{ on $\partial\Omega$.}\nonumber
\end{eqnarray}
Then by induction we conclude that $0\le\phi_{n}\le\phi_{n+1}$
since the function $g(s)$ is increasing. The same induction
argument shows that $\phi_n(x)\le\bar\phi(x)$ for all $n\ge 1$.
Therefore, the sequence $\phi_n$ converges to a limit
$\phi_\lambda$ which has to be a solution of (\ref{elliptic}) and
satisfy $0\le\phi_\lambda\le \bar\phi(x)$. As the sequence
$\phi_n$ does not depend on the choice of the super-solution
$\bar\phi$, the limit $\phi_\lambda$ is a minimal solution of
(\ref{elliptic}). $\Box$

{\bf Proof of Lemma \ref{prop-small-lambda}.} Observe that for
$\lambda>0$ sufficiently small the function
$\bar\tau(x)=2g(0)\lambda\tau(x)$ satisfies
\begin{eqnarray}\label{eq-bartau}
&&-\Delta\bar\tau+u\cdot\nabla\bar\tau\ge
\lambda g(\bar\tau)\hbox{ in $\Omega$,}\\
&&\bar\tau=0\hbox{ on $\partial\Omega$.}\nonumber
\end{eqnarray}
Here $\tau(x)$ is the solution of (\ref{eq-tau0}).
This is true provided that $2g(0)\ge g(2g(0)\lambda\tau)$. As the
function $g(s)$ is increasing, for this inequality to hold it suffices to require that
$2g(0)\ge g(2g(0)\lambda \theta_u)$. This condition is clearly
satisfied if $\lambda\le C/\theta_u$ with a constant $C$ which
depends only on the function $g(s)$. Now, existence of a minimal
solution to (\ref{elliptic}) follows from Lemma~\ref{lem-super}.
$\Box$

Recall that a solution $\phi_\lambda$ of (\ref{elliptic}) is
stable if the principal eigenvalue
$\kappa_1(\lambda,\phi_\lambda)$ of the linearized operator
\[
M_\lambda\psi=-\Delta\psi+u\cdot\nabla\psi-\lambda
g'(\phi_\lambda)\psi
\]
is positive.
\begin{lemma}\label{prop-stable-0}
Any minimal solution of (\ref{elliptic}) has
$\kappa_1(\lambda,\phi_\lambda)\ge 0$.
\end{lemma}
{\bf Proof.} Let $\phi_\lambda$ be a minimal solution of (\ref{elliptic})
and assume that the principal eigenvalue $\kappa_1(\lambda,\phi_\lambda)$
of the problem
\begin{eqnarray}\label{stable-eq1}
&&-\Delta\psi+u\cdot\nabla\psi-\lambda
g'(\phi_\lambda)\psi=\kappa_1(\lambda,\phi_\lambda)\psi,\\
&&\psi=0\hbox{ on $\partial\Omega$}\nonumber
\end{eqnarray}
is negative. Consider the function $\psi_\eps=\phi_\lambda-\eps\psi$, then
we have
\begin{eqnarray*}
&&-\Delta\psi_\eps+u\cdot\nabla\psi_\eps-\lambda
g(\psi_\eps)=\lambda g(\phi_\lambda)-\eps\lambda
g'(\phi_\lambda)\psi-\eps\kappa_1(\phi_\lambda)\psi-\lambda
g(\phi_\lambda-\eps\psi)\\
&&~~~~~~~~~~~~~~~~~~~~~~~~~~~~~~~~~~
=-\eps\kappa_1(\lambda,\phi_\lambda)\psi
-\frac{\eps^2g''(\xi)}{2}\psi^2\ge
0,
\end{eqnarray*}
provided that $\eps$ is sufficiently small and
$\kappa_1(\lambda,\phi_\lambda)< 0$. This contradicts minimality of
$\phi$. Therefore, we have $\kappa_1(\lambda,\phi_\lambda)\ge 0$ if $\phi_\lambda$ is a
minimal solution. $\Box$

\begin{lemma}\label{lem-kzero}
Assume that $\phi_\lambda$ is a solution of (\ref{elliptic}) such
that $\kappa_1(\lambda,\phi_\lambda)=0$. Then no classical
solutions of (\ref{elliptic}) with $\tilde\lambda>\lambda$ exists.
\end{lemma}
{\bf Proof.} Assume that $\tilde\lambda>\lambda$ and there exists
a function $\tilde\phi\ge 0$ such that
\begin{eqnarray*}
&&-\Delta\tilde\phi+u\cdot\nabla\tilde\phi=
\tilde\lambda g(\tilde\phi), \\
&&\tilde\phi=0\hbox{ on $\partial\Omega$.}
\end{eqnarray*}
Let us also denote by $\psi$ the positive eigenfunction of the
adjoint problem
\begin{eqnarray}\label{stable-eq11}
&&-\Delta\psi-\nabla\cdot(u \psi)-\lambda
g'(\phi_\lambda)\psi=0,\\
&&\psi=0\hbox{ on $\partial\Omega$.}\nonumber
\end{eqnarray}
Set $\eta=\phi_\lambda+\tau(\tilde\phi-\phi_\lambda)$ with
$\tau\in[0,1]$. Then convexity of $g$   implies that
\begin{eqnarray}\label{0-etaplus0}
&&-\Delta \eta+u\cdot\nabla\eta-\lambda g(\eta)=
-\Delta \eta+u\cdot\nabla\eta-
\lambda g(\phi_\lambda+\tau(\tilde\phi-\phi_\lambda))\\
&&~~~~~~~~~~~~~~~~~ \ge -\Delta \eta+u\cdot\nabla\eta -\lambda
(1-\tau)g(\phi_\lambda)-\lambda \tau
g(\tilde\phi)=(\tilde\lambda-\lambda)\tau g(\tilde\phi) \ge
0,\nonumber
\end{eqnarray}
for all $\tau\in[0,1]$. Moreover, we have equality in
(\ref{0-etaplus0}) when $\tau=0$.  Differentiating
(\ref{0-etaplus0}) with respect to $\tau$ at $\tau=0$ gives the
following inequality for $\zeta=\tilde\phi-\phi_\lambda$:
\begin{equation}\label{0-etaplus10}
-\Delta\zeta+u\cdot\nabla\zeta-\lambda g'(\phi_\lambda)\zeta\ge
(\tilde\lambda-\lambda) g(\tilde\phi)>0.
\end{equation}
Multiplying (\ref{0-etaplus10}) by the eigenfunction $\psi$ of
(\ref{stable-eq11}) and integrating we obtain
\[
0<\int\psi\left[-\Delta\zeta+u\cdot\nabla\zeta-\lambda
g'(\phi_\lambda)\zeta\right]=\int
\zeta\left[-\Delta\psi-\nabla\cdot(u\psi)-\lambda
g'(\phi_\lambda)\psi\right]=0.
\]
This contradiction shows that no solution of (\ref{elliptic}) for
$\tilde\lambda>\lambda$ may exist if
$\kappa_1(\lambda,\phi_\lambda)=0$. $\Box$

This also finishes the proof of Proposition~\ref{thm1}. The
critical threshold $\lambda^*(u)$ is the supremum of all $\lambda$
for which a stable solution of (\ref{elliptic}) exists. We
summarize the upper and lower bounds for $\lambda^*(u)$ in Lemmas
\ref{prop-large-lambda} and \ref{prop-small-lambda} as
\begin{equation}\label{upper-mu}
\frac{C}{\theta_u}\le\lambda^*(u)\le\frac{\mu_1(u)}{g'(0)}<+\infty.
\end{equation}
We will use these bounds in the sequel.

\subsection{A uniform bound away from $\lambda^*$}

Uniform $L^\infty$-bounds for the functions $\phi_{\lambda^*}$ at
$\lambda=\lambda^*$ are difficult to obtain and will be investigated
elsewhere~\cite{CBR}. However, we have the following universal
estimate for $\lambda<\lambda^*$ which will prove useful later.
\begin{prop}\label{prop-univ-bd}
For any $\delta>0$ there exists a constant $C(\delta)>0$ which depends
only on $\delta$ and nonlinearity $g(s)$ but not on the domain
$\Omega$ or the incompressible flow $u(x)$ so that the minimal
positive solution $\phi_\lambda(x)$ of (\ref{elliptic}) satisfies
$0\le\phi_\lambda(x)\le C(\delta)$ for all
$\lambda\in(0,(1-\delta)\lambda^*)$.
\end{prop}
{\bf Proof.} The proof is based on an idea from \cite{BCMR}. Fix
$\delta\in(0,1)$, let $\lambda_0<(1-\delta)\lambda^*$ and take any
$\lambda_1\in((1-\delta/3)\lambda^*,\lambda^*)$. We denote by $\phi_0$ and
$\phi_1$ the corresponding classical solutions of (\ref{elliptic})
with $\lambda=\lambda_0$ and $\lambda=\lambda_1$, respectively.

Following \cite{BCMR}, set
\[
h(s)=\int_0^s\frac{ds'}{g(s')}.
\]
It follows from the positivity of the function $g(s)$ and
(\ref{gint-infty}) that $h(s)$ is an increasing positive function with
$h(+\infty)<+\infty$. We now define the rescaled inverse function
\begin{equation}\label{def-Phi}
\Phi(s)=h^{-1}\left(\frac{\lambda_0}{\lambda_1}h(s)\right).
\end{equation}
Note that, since $\lambda_0<(1-\delta)\lambda^*$ and
$\lambda_1>(1-\delta/3)\lambda^*$ we have
\[
0\le\frac{\lambda_0}{\lambda_1}h(s)<
\frac{\lambda_0}{\lambda_1}h(+\infty)<\frac{1-\delta}{1-\delta/3}h(+\infty).
\]
Therefore, the function $\Phi(s)$ is well-defined for all $s\ge 0$,
and there exists a constant $K(\delta)$ which depends only on the parameter $\delta>0$
and the nonlinearity $g(s)$ so that $0\le\Phi(s)\le K(\delta)$ for all
$s\ge 0$ and all $\lambda\in(0,(1-\delta)\lambda^*)$.

In addition, as $g(s)\ge g(0)=1$, we have $\Phi(s)\le s$ and
\begin{equation}\label{basic-phiprime}
\Phi'(s)=\left[h'\left(h^{-1}\left(\frac{\lambda_0}{\lambda_1}h(s)\right)
\right)\right]^{-1}
 \frac{\lambda_0}{\lambda_1}h'(s)=
\frac{\lambda_0 g(\Phi(s))}{\lambda_1 g(s)}.
\end{equation}
Hence, as $g(s)$ is increasing and $\Phi(s)\le s$, the function
$\Phi(s)$ is increasing, with
\[
0<\Phi'(s)\le
\lambda_0/\lambda_1<1.
\]
Moreover, $\Phi$ is concave:
\begin{eqnarray*}
&&\!\!\!\!\Phi''(s)=
\frac{\lambda_0}{\lambda_1}
\frac{g'(\Phi(s))\Phi'(s)g(s)-g(\Phi(s))g'(s)}{g^2(s)}=
\frac{\lambda_0}{\lambda_1 g^2(s)}\left[
\frac{\lambda_0 g(\Phi(s))}{\lambda_1 g(s)}g'(\Phi(s))g(s)-g(\Phi(s))g'(s)
\right]\\
&&~~~~=
\frac{\lambda_0 g(\Phi(s)}{\lambda_1 g^2(s)}\left[
\frac{\lambda}{\lambda_1}g'(\Phi(s))-g'(s)\right]\le 0
\end{eqnarray*}
because $g(s)$ is convex, $\Phi(s)\le s$ and $0<\lambda<\lambda_1$.

Recall that $\phi_1$ is the minimal positive classical solution
to (\ref{elliptic}) with $\lambda=\lambda_1$ and set
$\bar\phi=\Phi(\phi_1)$. Using concavity of the function $\Phi(s)$ and
expression (\ref{basic-phiprime}) we observe that the function $\bar\phi$
satisfies the inequality
\begin{eqnarray*}
&&-\Delta\bar\phi+u\cdot\nabla\bar\phi=
-\Phi''(\phi_1)|\nabla\phi_1|^2
-\Phi'(\phi_1)\Delta\phi_1+\Phi'(\phi_1)(u\cdot\nabla\phi_1)
\\
&&=-\Phi''(\phi_1)|\nabla\phi_1|^2+\lambda_1\Phi'(\phi_1)g(\phi_1)
\ge
\lambda_1\Phi'(\phi_{1})g(\phi_{1})=
\lambda_0g(\bar\phi).
\end{eqnarray*}
Moreover, as $\Phi(0)=0$ the function $\bar\phi$ obeys the Dirichlet
boundary conditions $\bar\phi=0$ on $\partial\Omega$. Therefore, the
function $\bar\phi$ is a super-solution for (\ref{elliptic}) with
$\lambda=\lambda_0$. Employing the same iterative procedure as in the
proof of Proposition \ref{prop-small-lambda} we may then construct a
non-negative solution $\phi_0$ of (\ref{elliptic}) with
$\lambda=\lambda_0$ which is smaller than $\bar\phi(x)$. However, by
construction we have $0\le\bar\phi(x)\le K(\delta)$ and the conclusion
of Proposition \ref{prop-univ-bd} holds. $\Box$

\subsection{A uniform lower bound for $\lambda^*$}\label{sec:uniform}

We prove Theorem \ref{thm-lowerbd} in this section. Let $\phi(x)$ be the minimal
positive
solution of (\ref{elliptic}):
\begin{eqnarray}\label{lbd-1}
&&-\Delta\phi+u\cdot\nabla\phi=\lambda g(\phi)\hbox{ in $\Omega$,} \\
&&\phi=0\hbox{ on $\partial\Omega$.} \nonumber
\end{eqnarray}
According to (\ref{upper-mu}), in order to obtain a uniform lower bound
for the explosion threshold $\lambda^*$, it suffices to bound from
above $\theta_u$, the supremum of the exit time, defined by
(\ref{eq-tau0}) and (\ref{def-thetau}).  That is, it suffices to prove
that there exists a constant $M>0$ so that
\begin{equation}\label{lbd-unif}
\theta_u\le M,
\end{equation}
for all divergence free flows $u(x)$ in
$\Omega$. The constant $M$ should not depend on the flow $u(x)$.
This bound is an immediate consequence of Lemma~\ref{lem-lp-linfty}.

\subsubsection*{Proof of Lemma~\ref{lem-lp-linfty}}

We write $q(x)$, the solution of (\ref{lbd-lem}), as
\begin{equation}\label{lbd-tau-psi}
q(x)=\int_0^\infty\bar\psi(t,x)dt.
\end{equation}
The function $\bar\psi(t,x)$ satisfies the parabolic initial value problem
\begin{eqnarray}\label{lbd-3}
&&\bar\psi_t-\Delta\bar\psi+u\cdot\nabla\bar\psi=0\hbox{ in $\Omega$,} \\
&&\bar\psi(t,x)=0\hbox{ on $\partial\Omega$,} \nonumber\\
&&\bar\psi(0,x)=f(x)\hbox{ in $\Omega$.}\nonumber
\end{eqnarray}
We will now show that there exists a pair of constants $C>0$ and
$\alpha>0$ so that for any incompressible flow $u$ and any solution of
(\ref{lbd-3}) with initial data $f(x)$ we have a
uniform bound
\begin{equation}\label{lbd-4}
|\psi(t,x)|\le \frac{C e^{-\alpha t}}{t^r}\|f\|_{L^1(\Omega)},
\end{equation}
with any $r>d/2$.
The proof is as in \cite{CKRZ} with a slight modification, we present
the details for the convenience of the reader. First, multiplying
(\ref{lbd-3}) by $\psi$ and integrating by parts we obtain
\begin{equation}\label{lbd-5}
\frac 12\frac{d}{dt}\|\psi\|_2^2=-\|\nabla\psi\|_2^2.
\end{equation}
Using the Poincar\'e inequality in $\Omega$ we conclude that there
exists a constant $\alpha>0$ so that
\begin{equation}\label{lbd-51}
\|\psi(t_2)\|_2\le e^{-\alpha (t_2-t_1)}\|\psi(t_1)\|_2
\end{equation}
for any pair of times $t_2\ge t_1\ge 0$.
On the other hand, we have, using the Poincar\'e inequality again, for
all $1<p<2d/(d-2)$:
\[
\|\psi\|_p\le C\|\nabla\psi\|_2.
\]
Next, using the H\"older inequality, with $1/p+1/q=1$ we obtain:
\[
\|\psi\|_2^2=
\int|\psi|^2\le \left(\int|\psi|\right)^{1/p}
\left(\int|\psi|^{(2-1/p)q}\right)^{1/q}
\le C\|\psi\|_1^{1/p}\|\nabla\psi\|_2^{2-1/p},
\]
provided that
\[
\left(2-\frac{1}{p}\right)q=\left(2-\frac{1}{p}\right)\frac{p}{p-1}
=\frac{2p-1}{p-1}<\frac{2d}{d-2},
\]
or, equivalently, that $p>(d+2)/4$. Therefore, we have the following
Nash-type inequality in $\Omega$:
\[
\|\nabla\psi\|_2^2\ge C\frac{\|\psi\|_2^{4p/(2p-1)}}{\|\psi\|_1^{2/(2p-1)}}=
C\frac{\|\psi\|_2^{s+2}}{\|\psi\|_1^{s}},
\]
with $s=2/(2p-1)$. However, incompressibility of the flow, the Hopf lemma
and the boundary conditions imply that $\|\psi(t)\|_1\le \|f\|_1$.
It follows that
\[
\|\nabla\psi\|_2^2\ge
C\frac{\|\psi\|_2^{s+2}}{\|f\|_1^{s}}.
\]
Going back to (\ref{lbd-5}) we conclude that
\begin{equation}\label{lbd-6}
\frac{d}{dt}\|\psi\|_2^2=-2\|\nabla\psi\|_2^2\le
-C\frac{\|\psi\|_2^{s+2}}{\|f\|_1^{s}}.
\end{equation}
Therefore we have a bound
\[
\|\psi(t)\|_2\le \frac{C}{t^{1/s}}\|f\|_1.
\]
Combining this inequality with (\ref{lbd-51}) evaluated with $t_1=t$,
$t_2=2t_1$ we conclude that
\begin{equation}\label{lbd-7}
\|\psi(t)\|_2\le \frac{Ce^{-\alpha t}}{t^{1/s}}\|f\|_1,
\end{equation}
with $1/s>d/4$.

Consider now the solution operator ${\cal P}_t:~\psi_0\to\psi(t)$. We
have shown that
\[
\|{\cal P}_t\|_{L^1\to L^2}\le \frac{Ce^{-\alpha t}}{t^{1/s}}.
\]
The adjoint operator to ${\cal P}_t^*$ is simply the solution
operator corresponding to the (also incompressible) flow $(-u)$.
Therefore, we have the dual bound
\[
\|{\cal P}_t^*\|_{L^1\to L^2}\le \frac{Ce^{-\alpha t}}{t^{1/s}},
\]
which in turn implies that
\[
\|{\cal P}_t\|_{L^2\to L^\infty}\le \frac{Ce^{-\alpha t}}{t^{1/s}}.
\]
Putting these bounds together we obtain
\[
\|\psi(t)\|_\infty=\|{\cal P}_tf\|_\infty=
\|{\cal P}_{t/2}{\cal P}_{t/2}f\|_\infty
\le\|{\cal P}_{t/2}\|_{L^2\to L^\infty}
\|{\cal P}_{t/2}\|_{L^1\to L^2}\|f\|_1
\le \frac{Ce^{-\alpha t}}{t^{2/s}}\|f\|_1,
\]
which is (\ref{lbd-4}). The maximum principle also implies that we
have a trivial bound $\|\psi\|_{L^\infty}\le \|f\|_{L^\infty}$. Interpolating
between these two bounds we get the estimate
\begin{equation}\label{lin-lp-linfty}
|\psi(t,x)|\le \frac{C_\eps
e^{-\alpha_pt}}{t^{n/(2p)+\eps}}\|f\|_{L^p},
\end{equation}
for any $\eps>0$. Now, for any $p>n/2$ we may choose $\eps>0$
sufficiently small so that the kernel would be integrable at $t=0$, and
(\ref{lbd-tau-psi}) would imply that $\|q\|_{L^\infty}\le C\|f\|_{L^p}$ and
the constant $C>0$ is independent of the incompressible flow $u$.
This finishes the proof of Lemma~\ref{lem-lp-linfty} and hence that of
Theorem~\ref{thm-lowerbd}.  $\Box$

\subsubsection*{Explosion threshold in compressible flows}

As we have mentioned in the introduction, without the incompressibility
constraint the explosion threshold may be arbitrarily small. Indeed,
according to Proposition~\ref{prop-large-lambda} we have an upper bound
$\lambda^*(u)\le \mu_1(u)/g'(0)$. Therefore, to see that no uniform in the flow
lower bound for $\lambda^*$ in compressible flows exists, it suffices to
construct flows $u_n(x)$ such that the principal eigenvalue
$\mu_1(u_n)\to 0$, as $n\to+\infty$. Such example is provided by the
radial flows $u_n(x)=4n x$, say, in two-dimensions:
\begin{eqnarray}\label{comp-eig}
&&-\Delta\phi_n+4nx\cdot\nabla\phi_n=\mu_n\phi_n,\hbox{ $\phi_n>0$ in
$B(0,1)\subset\Rm^2$},\\
&&\phi_n=0\hbox{ on $|x|=1$}.\nonumber
\end{eqnarray}
Then $\mu_n\le C e^{-cn}\to 0$ as $n\to+\infty$ -- this
can be seen either from the general theory in \cite{Friedman,Kifer} or by an
explicit computation. Indeed, setting $\phi_n=e^{n|x|^2}\psi_n$
we obtain a self-adjoint problem for $\psi_n$:
\begin{eqnarray}
&&-\Delta\psi_n+4n^2 |x|^2\psi_n=(\mu_n+4n)\mu_n\psi_n,\hbox{ $\psi_n>0$ in
$B(0,1)\subset\Rm^2$},\\
&&\psi_n=0\hbox{ on $|x|=1$}.\nonumber
\end{eqnarray}
Hence, $\mu_n$ satisfies the variational principle
\[
\mu_n=-4n+\inf_{\psi\in H_0^1(B)}
\frac{\int|\nabla\psi|^2+4n^2\int|x|^2|\psi|^2}{\int|\psi|^2}.
\]
In addition, we have $\mu_n> 0$, as follows from the maximum principle
applied to (\ref{comp-eig})
For a test function of the form $\psi(x)=e^{-n|x|^2}q(x)$, where $0\le q(x)\le 1$,
$q(x)=1$ for $0\le|x|\le 1/2$ and $q(x)=0$ for $|x|\ge 3/4$ we obtain
by a straightforward computation
\[
\frac{\int|\nabla\psi|^2+4n^2\int|x|^2|\psi|^2}{\int|\psi|^2}
=8 n^2\frac{  \int|x|^2|\psi|^2}{\int|\psi|^2}+O(e^{-cn})= 8 n^2
\frac{  \int_{\mathbb{R}^2}|x|^2e^{-2n|x|^2}}{\int_{\mathbb{R}^2}e^{-2n|x|^2}}+O(e^{-cn})
=4n+O(e^{-cn}),
\]
with $c>0$.
Therefore, $\mu_n\to 0$ as $n\to +\infty$ and hence $\lambda_n^*\to 0$ as well.

\section{The strong flow asymptotics}

In this section we consider the elliptic problem (\ref{elliptic}) when
the advecting flow is strong. Accordingly, we introduce a large parameter $A\gg
1$ and re-write (\ref{elliptic}) as
\begin{eqnarray}\label{strong-elliptic}
&&-\Delta\phi_\lambda+Au\cdot\nabla\phi_\lambda=\lambda g(\phi_\lambda)
\hbox{ in $\Omega$},\\
&&\phi_\lambda=0\hbox{ on $\partial\Omega$.}\nonumber
\end{eqnarray}
We are interested in the behavior of the solution $\phi(x)$ of
(\ref{strong-elliptic}) for large $A$, as well as in the
dependence of the explosion threshold $\lambda^*$ on the amplitude
$A$. With a slight abuse of notation we will denote here by
$\lambda^*(A)$ the explosion threshold of the problem
(\ref{strong-elliptic}).

\subsection{Equidistribution on the flow streamlines}

Our first result shows that, when the flow is strong, solution becomes
nearly constant on the flow streamlines, at least in an average sense
and for $\lambda$ away from $\lambda^*(A)$. This is a common phenomenon
in diffusion-advection problems: a strong flow induces stratification.
\begin{prop}\label{thm7}
Assume in that $u\cdot n=0$ on the boundary $\partial\Omega$.
Then the solution $\phi_\lambda$ of (\ref{elliptic}) is nearly
constant on the streamlines of $u$ for sufficiently large $A$ in the
sense that for any $\delta>0$ there exists $C(\delta)$ so that
\[
\int|u\cdot\nabla\phi_\lambda|^2\le\frac{C(\delta)\lambda}{A},
\]
for all $\lambda\le(1-\delta)\lambda^*(A)$.
\end{prop}
{\bf Proof.} First, we multiply (\ref{strong-elliptic}) by
$u\cdot\nabla\phi_\lambda$ and integrate over $\Omega$. The emerging integrals
over the boundary vanish since $\phi_\lambda=u\cdot n=0$ on
$\partial\Omega$, and we obtain the following estimate:
\begin{eqnarray*}
&&\int|u\cdot\nabla\phi_\lambda|^2=\frac{1}{A}\int(u\cdot\nabla\phi_\lambda)
[\Delta \phi_\lambda+\lambda g(\phi_\lambda)]=
\frac{1}{A}\int(u\cdot\nabla\phi_\lambda)\Delta \phi_\lambda
=
-\frac{1}{A}\int_\Omega\pdr{\phi_\lambda}{x_k}\pdr{}{x_k}
\left[u_j\pdr{\phi_\lambda}{x_j}\right]\\
&&=
-\frac{1}{A}\int_\Omega\pdr{\phi_\lambda}{x_k}
\pdr{u_j}{x_k}\pdr{\phi_\lambda}{x_j}-
\frac{1}{A}\int_\Omega\pdr{\phi_\lambda}{x_k}
u_j\frac{\partial^2\phi_\lambda}{\partial x_j\partial x_k}=
-\frac{1}{A}\int_\Omega\pdr{\phi_\lambda}{x_k}
\pdr{u_j}{x_k}\pdr{\phi_\lambda}{x_j}
\le \frac{C}{A}\int|\nabla\phi_\lambda|^2.
\end{eqnarray*}
This means that the variation along the streamlines is smaller than
across.  Now, we have to bound the $L^2$-norm of
$\nabla\phi_\lambda$. However, multiplying (\ref{strong-elliptic}) by
$\phi_\lambda$ and integrating by parts we see that
\[
\int|\nabla\phi_\lambda|^2=\lambda\int g(\phi_\lambda)\phi_\lambda.
\]
Moreover, Proposition \ref{prop-univ-bd} shows that $\phi_\lambda$
satisfies a uniform bound $0\le\phi_\lambda\le C(\delta)$ as long as
$\lambda\in(0,(1-\delta)\lambda^*)$.  Therefore, for such $\lambda$ we
know that
\begin{equation}\label{strong-gradient}
\int|\nabla\phi|^2\le C\lambda,
\end{equation}
for all $A>0$. It follows that
\[
\int|u\cdot\nabla\phi|^2\le \frac{C\lambda}{A}.
\]
In that sense solution becomes uniform over the streamlines. $\Box$

As a consequence, for each fixed
$\lambda<\limsup_{A\to\infty}\lambda^*(A)$ we know that
\[
\int|u\cdot\nabla\phi_\lambda|^2\to 0, \hbox{ as $A\to+\infty$.}
\]
We will improve this statement for two-dimensional cellular flows in
Section \ref{sec:cellflow}.

An interesting by-product of the estimate (\ref{strong-gradient}) is
that there are no boundary or internal layers in this problem unlike
in the problems with boundary forcing in cellular
flows~\cite{Childress,FP,Heinze,Koralov,NPR,Rosenbluth,Shraiman,Soward}.
The reason is that $\phi_\lambda$ is
set to be constant on the boundary and the normal component of the
flow vanishes on the boundary -- hence, there is no "conflict" between
the boundary data and uniformization along the streamlines of the
flow.

\subsection{The critical parameter in the limit of a strong flow}

We prove here Theorem \ref{thm5}.  We recall that the assumption that
there is no $H_0^1(\Omega)$ first integral $H(x)$ such that
$u\cdot\nabla H=0$ almost everywhere is equivalent to the fact that
the principal Dirichlet eigenvalue $\mu_1(A)$ of the operator
$-\Delta+A\cdot\nabla$ on $\Omega$ tends to infinity as $A\to+\infty$
\cite{BHN}.

First, assume that $\mu_1(A)$ is bounded as
$A\to+\infty$. Then the upper bound in (\ref{upper-mu}) for
$\lambda^*(A)$ implies that $\limsup_{A\to+\infty}\lambda^*(A)$ is
finite as well.

Next, we show that if $u$ has no first integral in $H_0^1(\Omega)$
then $\lambda^*(A)\to+\infty$ as $A\to+\infty$. The proof is based on
the lower bound for $\lambda^*(A)$ in (\ref{upper-mu}).  The next
lemma is contained in \cite{Kifer} -- we provide a proof in the spirit
of \cite{CKRZ}.
\begin{lemma}\label{lem-tau-A}
We have $\theta_A\to 0$ as $A\to+\infty$ if $u$ has no first integrals in
$H_0^1(\Omega)$.
\end{lemma}
{\bf Proof.} As in the proof of  Theorem \ref{thm-lowerbd}
we represent the function $\tau_A(x)$ using the Duhamel formula as
\[
\tau_A(x)=\int_0^\infty \psi(t,x)dt,
\]
with the function $\psi(t,x)$ that solves the parabolic problem
\begin{eqnarray}\label{parabolic11}
&&\psi_t-\Delta\psi+Au\cdot\nabla\psi=0\hbox{ in $\Omega$},\\
&&\psi(0,x)=1,\nonumber\\
&&\psi=0\hbox{ on $\partial\Omega$.}\nonumber
\end{eqnarray}
It follows from the Theorem 5.3 of \cite{CKRZ} that if $u$ has no first
integrals in $H_0^1(\Omega)$ then for any $t_0>0$ we can find a flow
amplitude $A_0(t_0)>0$ so that $\|\psi(nt_0)\|_{L^\infty(\Omega)}\le
2^{-n}|\Omega|$ for all $A\ge A_0$. Therefore, we have an upper bound
\[
|\tau_A(x)|\le
t_0|\Omega|\sum_{n=0}^\infty\left(\frac{1}{2}\right)^n={2t_0|\Omega|}
\hbox{ for $A\ge A_0(t_0)$}.
\]
We conclude that $\|\tau_A\|_{L^\infty(\Omega)}\to 0$ as
$A\to+\infty$. $\Box$

This also finishes the proof of Theorem \ref{thm5}. $\Box$

\subsection{The explosion problem in a two-dimensional cellular flow}
\label{sec:cellflow}

\begin{figure}[ht!]
\centerline{\includegraphics[height=9.0cm]{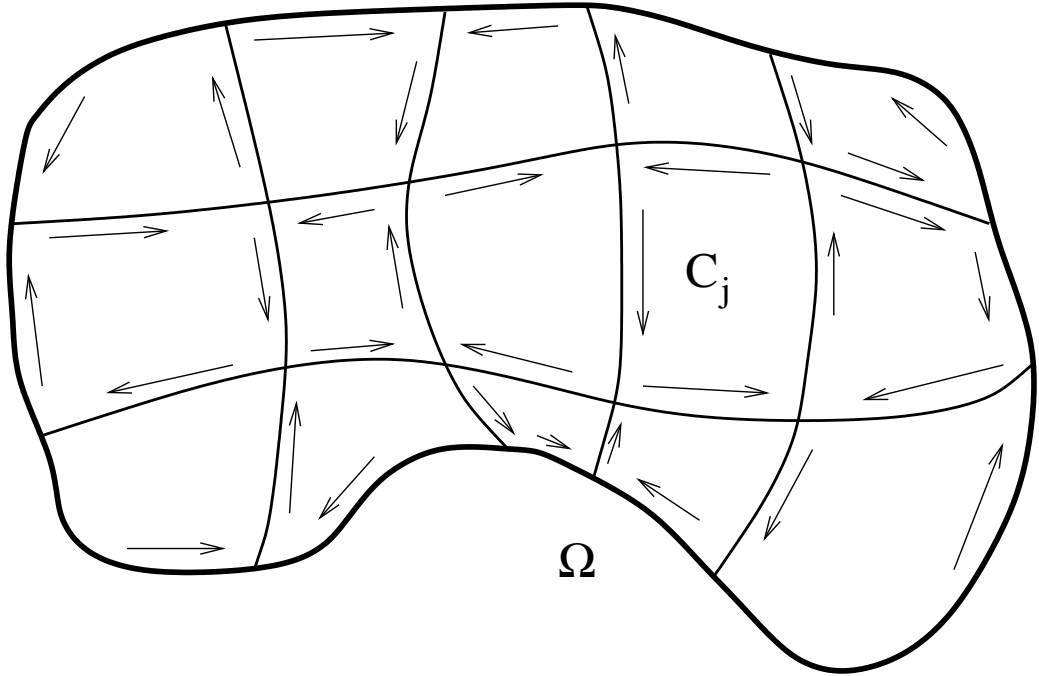}}
\caption{ A schematic description of a cellular flow.  }
\label{fig-omega-cell}
\end{figure}

We now consider the explosion problem
\begin{eqnarray}\label{cf-fr-1}
&&-\Delta\phi_A+Au\cdot\nabla\phi_A=\lambda f(\phi_A)
\hbox{ in $\Omega$,} \\
&&\phi_A=0\hbox{ on $\partial\Omega$,} \nonumber
\end{eqnarray}
in a two-dimensional simply connected domain $\Omega$. The flow has
the form $u=(\Psi_y,-\Psi_x)$ with a stream-function $\Psi(x,y)$ which
we assume to be sufficiently smooth. We assume that the boundary of
the domain $\Omega$ is a level set $\{\Psi=0\}$ which may contain
finitely many saddle critical points of the function $\Psi$~-- thus,
the boundary is a union of streamlines of the flow away from the
critical points. We also assume that inside $\Omega$ the flow has a
cellular structure: the saddles of $\Psi$ are all non-degenerate and
are connected by the flow separatrices which divide $\Omega$ into a
finite number of invariant regions, called the flow cells, that we
will denote by ${\cal C}_j$. The stream-function $\Psi(x,y)$ has only
one critical point $(x_0,y_0)$ inside each of ${\cal C}_j$, which is a
non-degenerate maximum or minimum.  A prototype example of such flow
has the stream-function $\Psi(x,y)=\sin \pi x\sin \pi y$ -- its cells
are squares $[n,n+1]\times[m,m+1]$ with integer $m$ and $n$, and the
domain $\Omega$ is a finite union of such squares. A more general flow
of such type is depicted in Figure~\ref{fig-omega-cell}.

\subsubsection*{The Freidlin problem}

The strong flow asymptotics for parabolic reaction-diffusion equations
for two-dimensional flows with Morse class stream-functions has been
considered in \cite{Freidlin}. This class does not include the
cellular flows under our consideration as we allow the stream-function to
have several saddles on the level set $\{\Psi=0\}$. Nevertheless, the
limit problem of \cite{Freidlin} is crucial in the explosion problem
in a cellular flow.  The limit problem in \cite{Freidlin}
was formulated as a system of reaction-diffusion equations on the Reeb
graph of the function $\Psi$. We recall and re-derive these
results below in the context of the explosion problem in the one-cell setting, as that is
what we will need below. The single cell is also the situation
addressed numerically in \cite{BKJS}.  For a one-cell flow the Reeb
graph is simply an interval $[0,H_0]$, where $H_0$ is the value of
$\Psi$ at the critical point inside the cell ${\cal C}$ that we assume
to be a maximum, and $\{\Psi=0\}$ is the boundary of the cell.  We are
interested in the behavior of solutions and of the explosion threshold
in the limit of a large flow amplitude.

The effective Freidlin problem on the interval $0\le h\le H_0$ is to
find a function $\bar\phi(h)$ satisfying
\begin{eqnarray}\label{cf-fr-3}
&&-\frac{1}{T(h)}\frac{d}{dh}\left(p(h)\frac{d\bar\phi}{dh}\right)=
\lambda g(\bar\phi),\\
&&\hbox{ $\bar\phi(0)=0$, $\bar\phi'(h)$ is bounded for $0\le h\le H_0$,}
\nonumber
\end{eqnarray}
with the coefficients
\begin{equation}\label{cf-fr-one-coeff}
T(h)=\oint_{\Psi(x,y)=h}\frac{dl}{|\nabla\Psi|},
~~p(h)=\oint_{\Psi(x,y)=h}{|\nabla\Psi|}dl.
\end{equation}
Under our assumptions on the stream-function,
the average turnover time $T(h)$ is bounded from above and
below away from zero:
\begin{equation}\label{cf-fr-T-bds}
0<T_0\le T(h)\le T_1|\ln h|.
\end{equation}
The uniform bound from below by $T_0$ in (\ref{cf-fr-T-bds}) comes
from the fact that the maximum of $\Psi(x,y)$ is a non-degenerate
critical point. The term $O(|\ln h|)$ for small $h>0$ appears in
(\ref{cf-fr-T-bds}) because the boundary may contain non-degenerate
saddle points of $\Psi$ so that the turnover time blows up as
$h\downarrow 0$.  The coefficient $p(h)$ is positive for $h>0$ and
behaves as $p(h)\sim C(H_0-h)$ close to $h=H_0$.  In particular we
have $p(H_0)=0$ (diffusivity vanishes at this point), while the drift
satisfies
\[
p'(h)=\oint_{\Psi(x,y)=h}
\frac{\Delta\psi}{|\nabla\Psi|}dl\le-\alpha_0,~~\hbox{with $\alpha_0>0$},
\]
for $h$ near $H_0$, and points away from $h=H_0$. Therefore, the
end-point $h=H_0$ is inaccessible for the diffusion process
corresponding to the left side of (\ref{cf-fr-3}), and one does not need to
prescribe the boundary condition at $h=H_0$ in order for
(\ref{cf-fr-3}) to be well-posed. The following proposition defines
the explosion threshold for the effective problem.
\begin{prop}\label{cf-fr-prop1}
There exists $\bar\lambda^*>0$ so that a positive solution of the
effective problem (\ref{cf-fr-3}) exists for all
$0\le\lambda<\bar\lambda^*$ and there is no positive solution of
(\ref{cf-fr-3}) for $\lambda>\bar\lambda^*$.
\end{prop}
{\bf Proof.} The proof follows the same steps as in
Section~\ref{sec:basic} -- the only required modification is due to
the degeneracies at $h=0$ and $h=H_0$. This can be addressed using the
general theory in \cite{Feller} and \cite{Mandl} but in the present
case the boundary value problem with a prescribed right side
\begin{eqnarray}\label{cf-fr-312}
&&-\frac{1}{T(h)}\frac{d}{dh}\left(p(h)\frac{d\psi}{dh}\right)=f(h)\\
&&\hbox{ $\psi(0)=0$, $\psi'(h)$ is bounded for $0\le h\le H_0$,}
\nonumber
\end{eqnarray}
has an explicit unique solution
\[
\psi(h)=\int_0^{h}\frac{1}{p(s)}\left(\int_s^{H_0}f(\xi)T(\xi)d\xi\right)ds
=\int_0^{H_0} f(\xi)T(\xi)P(\min{(h,\xi)})d\xi,
\]
where
\[
P(\xi)=\int_{0}^{\xi}\frac{ds}{p(s)}.
\]
In particular, we have $|P(\xi)|\le C|\ln(H_0-\xi)|$, so that
\[
|\psi(h)|\le \int_0^{H_0} |f(\xi)|T(\xi)P(\min{(h,\xi)})d\xi\le
C\|f\|_{\infty}\int_0^{H_0}|\ln \xi||\ln(H_0-\xi)|d\xi\le C\|f\|_\infty,
\]
and we also have
\begin{eqnarray*}
&&|\psi'(h)|\le\frac{1}{p(h)}\int_h^{H_0}|f(\xi)|T(\xi)d\xi\le
\frac{C\|f\|_\infty}{H_0-h}
\int_h^{H_0}|\ln \xi |d\xi \le {C\|f\|_\infty},
\end{eqnarray*}
so that
$\|\psi\|_{W^{1,\infty}}\le
C\|f\|_{L^\infty}$. Therefore, the mapping $f(h)\to \psi(h)$ is a
compact map  on $C[0,H_0]$ and the Krein-Rutman theory applies to
the operator in the left side of (\ref{cf-fr-312}). We may then
repeat the proof of Proposition~\ref{thm1} essentially verbatim
and conclude that the critical threshold $\bar\lambda$ for
(\ref{cf-fr-3}) exists. $\Box$

\subsubsection*{The explosion threshold
for strong cellular flows: the main result}

\begin{figure}[ht!]
\centerline{\includegraphics[height=9.0cm]{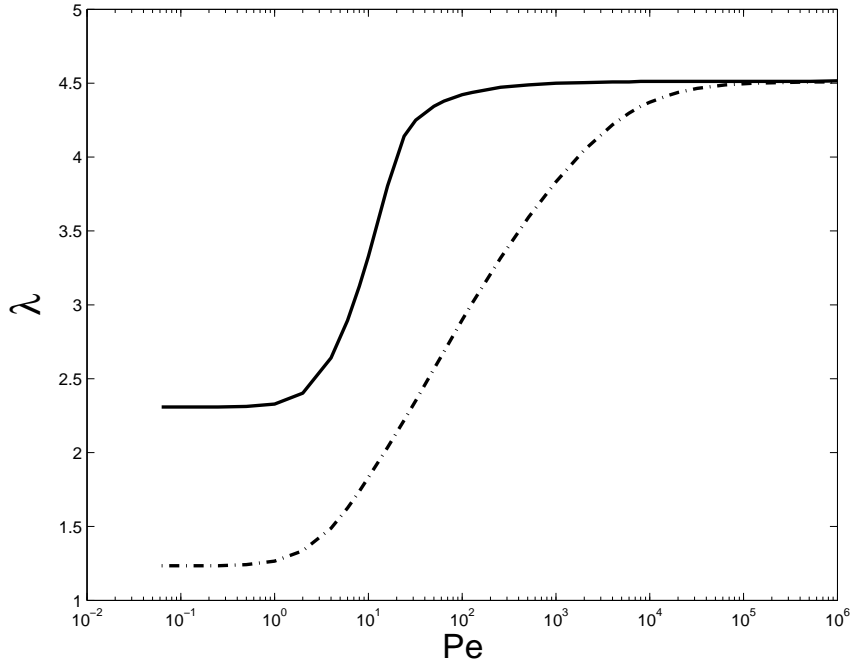}}
\caption{The value of $\lambda^*$ for various values of the P\'eclet number
$Pe=A$ for the cellular flow with the stream-function $\psi(x,y)= \sin
(2\pi ((2x/3+1)^3/8+1)) \sin(2\pi ((y+1)^2/4+1))$ on the the domain
$[0,2\pi]\times[0,2\pi]$ with four cells.  Dashed-dotted line:
$\lambda^*$ for the whole domain $[0,2\pi]\times[0,2\pi]$, solid line
-- the minimum of $\lambda^*$ of the four individual cells.  }
\label{cf-fig-lambda-many-cells}
\end{figure}

The main result of this section is the following theorem. We assume
that the flow has a cellular structure and satisfies the assumptions
outlined at the beginning of this section. Then for each cell ${\cal
C}_j$ one may formulate the corresponding one-cell Freidlin problem
(\ref{cf-fr-3}) for a function $\bar\phi_j$, posed now on an interval
$[0,H_j]$, where the outer boundary of ${\cal C}_j$ is the level set
$\{\Psi=0\}$ and $H_j$ is the value of the function $\Psi$ at the
(unique) extremal point inside ${\cal C}_j$. For the Freidlin problem
the Dirichlet boundary condition $\bar\phi(0)=0$ is prescribed at
$h=0$, and the derivative $\bar\phi_j'(H_j)$ is imposed to be
bounded. This defines the explosion threshold $\bar\lambda_j^*$ for
each cell ${\cal C}_j$. The following theorem shows that in the limit
of a large flow the explosion threshold for the whole domain $\Omega$
approaches the Freidlin explosion threshold for the "largest" cell
${\cal C}_j$.
\begin{thm}\label{cf-fr-thm}
Let $\lambda^*(A)$ be the explosion threshold for (\ref{cf-fr-1}) and
$\bar\lambda_j^*$ be the threshold for the aforementioned effective
problem (\ref{cf-fr-3}) posed in the cell ${\cal C}_j$ of the domain
$\Omega$.  Then we have
\[
\lim_{A\to\infty}\lambda^*(A)=\min_j\bar\lambda_j^*.
\]
\end{thm}
A numerical illustration of the main result of Theorem~\ref{cf-fr-thm}
is depicted in Figure~\ref{cf-fig-lambda-many-cells}.

The proof of Theorem~\ref{cf-fr-thm} proceeds in several steps. First,
we prove a stratification lemma for solutions of forced
advection-diffusion problems in cellular flows. It shows that solution of the Dirichlet
problem is small not only on the outer boundary but also on the whole skeleton of separatrices
and cells "do not talk to each other".
The second step is to
establish the result of Theorem~\ref{cf-fr-thm} for domains consisting
of one cell where one just has to show that for one cell the
explosion threshold converges in the strong flow limit to that of the
Freidlin problem. The last step is to generalize this result to a
domain consisting of finitely many cells.

\subsection{Cellular flows: a stratification lemma}

We begin the proof of Theorem~\ref{cf-fr-thm} with the following
lemma, of an independent interest.  Let ${\cal D}_0$ be the union of
all cell boundaries (separatrices of the flow) inside $\Omega$
including the outside boundary $\partial\Omega$.  First, we show that
solutions of a linear problem with the homogeneous Dirichlet data on the outer
boundary are small on ${\cal D}_0$.
\begin{lemma}\label{cf-lem-exit}
Let $\psi_A(x)$ be the exit time from $\Omega$, solution of
\begin{eqnarray}\label{cf-out-14}
&&-\Delta\psi_A+Au\cdot\nabla\psi_A=1
\hbox{ in $\Omega$},\\
&&\psi_A= 0\hbox{ on $\partial\Omega$.}\nonumber
\end{eqnarray}
For any $\delta>0$ there exists $A_0>0$ so that for all $A>A_0$
we have $0\le\psi_A(x)\le\delta$ for all $x\in{\cal D}_0$.
\end{lemma}
Intuitively, this lemma says that once a diffusive particle obeying
an SDE
\[
dX_t=Au(X_t)dt+\sqrt{2}dW_t,
\]
comes close to the skeleton of separatrices, somewhere inside
$\Omega$, it exits the domain $\Omega$ after a short time.  This is
the phenomenon behind the effective
diffusivity
\cite{Childress,FP,Heinze,Koralov,NPR,Rosenbluth,Shraiman,Soward},
and front speed and principle eigenvalue
enhancement~\cite{ABP,KR,NR,RZ,Zlatos} in cellular flows.

{\bf Proof.} The proof is in two steps.  First, we show that for any
$\delta>0$ there exists  a small $h>0$ and a large $A_0>0$ so that for any $A>A_0$ we
can find $h_j(A)$ with $h\le|h_j|\le 2h$, such that
$0\le\psi_A(x)\le\delta/2$ for all $x$ on the streamline $\{x\in{\cal
C}_j:~\Psi(x)=h_j(A)\}$ inside the cell ${\cal C}_j$. Let $D_j(A)$ be
the interior of those streamlines. In the second step we consider the
water-pipe domain $P_A=\Omega\setminus\left(\bigcup_j D_j(A)\right)$,
that is, a narrow tube around the skeleton of separatrices.  Using the
fact that $h_j$ is small, we apply the maximum principle for narrow
domains to conclude that the function $\psi_A$ is smaller than
$\delta$ in all of $P_A$ and not only on its boundary. As a
consequence, $\psi_A$ is small also on the skeleton ${\cal D}_0$.

{\bf Step 1.} We have a uniform $L^2$-bound for the gradient:
\begin{equation}\label{cf-out-grad-bd}
\int_\Omega|\nabla\psi_A(x)|^2dx\le C,
\end{equation}
with the constant $C$ independent of $A>0$, which follows from
multiplying (\ref{cf-out-14}) by $\psi_A(x)$ and integrating by parts,
together with the uniform $L^\infty$-bound for $\psi_A(x)$:
$\|\psi_A\|_{L^\infty}\le C$, which follows from
Lemma~\ref{lem-lp-linfty}.

Now, take $s\in(0,h/4)$ and let $F_j(h,s)$ be
the domain between the two streamlines $\{\Psi(x)=5h/4-s\}$ and
$\{\Psi(x)=7h/4+s\}$ inside the cell ${\cal C}_j$. We multiply
(\ref{cf-out-14}) by $(u\cdot\nabla\psi_A)$ and integrate over $F_j$:
\begin{eqnarray*}
&&\int_{F_j}|u\cdot\nabla\psi_A|^2=
\frac{1}{A}\int_{F_j}(u\cdot\nabla\psi_A)\Delta\psi_A dx
=\frac{1}{A}\int_{F_j}\sum_{m,k}u_m
\pdr{\psi_A}{x_m}\frac{\partial^2\psi_A}{\partial x_k^2}dx\\
&&=
\frac{1}{A}\int_{\partial F_j}(u\cdot\nabla\psi_A)(n\cdot\nabla\psi_A)dl-
\frac{1}{A}\int_{F_j}\pdr{u_m}{x_k}\pdr{\psi_A}{x_m}\pdr{\psi_A}{x_k}dx\le
\frac{C}{A}\int_{\partial F_j}|\nabla\psi_A(x)|^2dx+\frac{C}{A}.
\end{eqnarray*}
We used (\ref{cf-out-grad-bd}) in the last step.
Averaging this estimate in $s\in(0,h/4)$ we conclude that
\[
\int_{\bar F_j}|u\cdot\nabla\psi_A|^2dx\le \frac{C(h)}{A},
\]
where $\bar F_j=F_j(h,0)$ is the domain between the streamlines
$\{\Psi(x)=5h/4\}$ and $\{\Psi(x)=7h/4\}$ inside the cell ${\cal
C}_j$. The constant $C(h)$ may blow-up as $h\downarrow 0$ but that is
not important at the moment.

It follows that there exists a value
$h_j(A)\in(5h/4,7h/4)$ so that along the streamline
$L_j(A)=\{x\in{\cal C}_j:~\Psi(x)=h_j(A)\}$ we have
\[
\oint_{L_j(A)}|u\cdot\nabla\psi_A|^2dl\le \frac{C(h)}{A},
\]
with a new constant $C(h)$. Therefore, the oscillation of $\psi_A$
along $L_j(A)$ is small:
\[
\hbox{osc}_{L_j(A)}\psi_A(x)\le \frac{C(h)}{\sqrt{A}}.
\]
Hence, $\psi_A(x)$ is close to a constant $M_j(A)$ on the
streamline $L_j(A)$ when $A$ is sufficiently large. As a consequence
of the gradient bound (\ref{cf-out-grad-bd}), we have
$|M_j(A)-M_m(A)|\le C\sqrt{h}$ if the cells ${\cal C}_j$ and ${\cal
C}_m$ have a common piece of the boundary.  As the outer boundary
$\partial\Omega$, where $\psi_A(x)=0$, is also part of some cell
boundaries, it follows that $0\le M_j(A)\le C\sqrt{h}$ for all cells
${\cal C}_j$.  Therefore, we have
\[
0\le\psi_A(x)\le C\sqrt{h}+\frac{C(h)}{\sqrt{A}}<\frac{\delta}{2},
\hbox{ for $x\in L_j(A)$}
\]
if
$h\in(0,h_0)$ is sufficiently small and $A>A_0$ is large enough.

{\bf Step 2.} Now, we look at the water-pipe ${P}_A$ and show that
solution is below $\delta$ everywhere in $P_A$. The function $\psi_A$
inside $P_A$ satisfies $0\le\psi_A(x)\le\delta/2+r_A$, where $r_A$ is
the exit time from the slightly larger domain
$Q_{2h}=\Omega\setminus\left(\bigcup_j D_j(2h)\right)$:
\begin{eqnarray}\label{cf-exit-Peps}
&&-\Delta r_A+Au\cdot\nabla r_A=1
\hbox{ in $Q_{2h}$},\\
&&r_A= 0\hbox{ on $\partial Q_{2h}$.}\nonumber
\end{eqnarray}
Now, as in the proof of Lemma~\ref{lem-lp-linfty} we conclude that
there exists a constant $C(h)$ such that
$\|r_A\|_{L^\infty(Q_{2h})}\le C(h)$ for all $A>0$.  The same proof
shows that $C(h)\to 0$ as $h\to 0$ -- this happens because the
principal eigenvalue of the Dirichlet Laplacian in $Q_{2h}$ tends to
infinity as $h\to 0$ while the constants $K_p(h)$ in the Poincar\'e
inequality $\|\psi\|_{L^p(Q_{2h})}\le
K_p(h)\|\nabla\psi\|_{L^2(Q_{2h})}$, $1<p<\infty$, satisfy $K_p\to 0$
as $h\to 0$. It follows that if we take $h>0$ sufficiently small
(independent of $A$) then $0\le r_A(x)\le\delta/2$ for all $A>0$.
Therefore, we have $0\le\psi_A(x)\le\delta$ for all $A>A_0$ in $P_A$
and, in particular, $0\le\psi_A(x)\le\delta$ on ${\cal D}_0$. The
proof of Lemma~\ref{cf-lem-exit} is now complete. $\Box$

\subsection{Explosion problem in a one-cell domain}\label{sec:freidlin}

We now consider the explosion problem
\begin{eqnarray}\label{fr-1}
&&-\Delta\phi_A+Au\cdot\nabla\phi_A=\lambda f(\phi_A)
\hbox{ in $\Omega$,} \\
&&\phi_A=0\hbox{ on $\partial\Omega$,} \nonumber
\end{eqnarray}
in a domain $\Omega$ which consists of just one flow cell.  Without
loss of generality we assume that the single critical point
$(x_0,y_0)\in\Omega$ of $\Psi$ inside $\Omega$ is a maximum and set
$H_0=\Psi(x_0,y_0)$. Let us now formulate the version of
Theorem~\ref{cf-fr-thm} for a one-cell domain.
\begin{prop}\label{cf-fr-thm-onecell}
Let $\Omega$ be a one-cell domain and let $\lambda^*(A)$ be the
explosion threshold for (\ref{fr-1}) and $\bar\lambda^*$ be the
threshold for the aforementioned effective problem (\ref{cf-fr-3})
posed on $[0,H_0]$.  Then we have
\[
\lim_{A\to\infty}\lambda^*(A)=\bar\lambda^*.
\]
\end{prop}

The proof of Proposition~\ref{cf-fr-thm-onecell} is in two
steps. First, passing from the problem on the cell to the Freidlin
problem we show that the Fredlin threshold $\bar\lambda^*$ is not
smaller than than $\limsup_{A\to+\infty}\lambda^*(A)$. Next, we
establish the opposite inequality by starting with a solution to the
Freidlin problem and constructing a super-solution for (\ref{fr-1}).
The second step is quite straightforward in the case when the boundary
of $\Omega$ contains no saddles of the flow $u$ but is somewhat more technical
if $\partial\Omega$ contains such fixed points.

\subsubsection*{Passage from the cell to the Freidlin problem}

We first prove that
\begin{equation}\label{fr-lambda-lambda-bar}
\limsup_{A\to\infty}\lambda^*(A)\le \bar\lambda^*.
\end{equation}
Assume that
\begin{equation}\label{fr-lambda-lambda-bar1}
\lambda<\limsup_{A\to\infty}\lambda^*(A).
\end{equation}
We will show that then $\lambda<\bar\lambda^*$ by constructing a
solution to the Freidlin problem (\ref{cf-fr-3}) as the limit
of a sequence of problems on $\Omega$.  It follows from
(\ref{fr-lambda-lambda-bar1}) that there exists $\delta>0$ and a
sequence $A_n\to +\infty$ such that
$\lambda<(1-\delta)\lambda^*(A_n)$. Therefore, as a consequence of
Proposition~\ref{prop-univ-bd},
the minimal positive solutions of
\begin{eqnarray}\label{fr-1-n}
&&-\Delta\phi_n+A_nu\cdot\nabla\phi_n=
\lambda f(\phi_n)\hbox{ in $\Omega$,} \\
&&\phi_n=0\hbox{ on $\partial\Omega$} \nonumber
\end{eqnarray}
are uniformly bounded in $L^\infty(\Omega)\cap H_0^1(\Omega)$:
\begin{equation}\label{fr-uniform}
0\le\phi_n\le C,~~\int|\nabla\phi_n|^2dx\le C,
\end{equation}
with the constant $C>0$ independent of $n$. Hence, the sequence
$\phi_n$ converges weakly in $H^1(\Omega)$ (after extracting a
subsequence) and strongly in $L^2(\Omega)$ to a function
$\bar\phi$. As the functions $\phi_n$ are uniformly bounded and
$g(\phi)$ is smooth, the sequence $g(\phi_n)$ converges to
$g(\bar\phi)$.

We claim that $\bar\phi$ depends only on the variable $h=\Psi(x,y)$ and
satisfies the effective Freidlin problem (\ref{cf-fr-3}). The first claim
follows after we divide (\ref{fr-1-n}) by $A_n$ and let $n\to
+\infty$. This leads to
\begin{equation}\label{cf-unablaphi}
u\cdot\nabla\bar\phi=0
\end{equation}
in the sense of distributions.  It is convenient now to introduce the
curvilinear coordinates $(h,\theta)$. The coordinates are chosen
so that $h(x,y)=\Psi(x,y)$, that is, the streamlines of the flow
are $\{h=\hbox{const}\}$, and the level lines of the coordinate
$\theta=\Theta(x,y)$ are orthogonal to the flow lines:
$\nabla\Theta\cdot\nabla\Psi=0$. We normalize $\theta$ so that
$0\le\theta\le 2\pi$ and the boundary $\partial\Omega$ is a level
set: $\partial\Omega=\{h=0\}$.
Then (\ref{cf-unablaphi})  implies that $\bar\phi$
depends only on the variable $h$. The $L^\infty$-bound in
(\ref{fr-uniform}) implies that $0\le\bar\phi(h)\le C$. In
addition, we have
\begin{eqnarray*}
\int|\nabla_x\bar\phi|^2dx=\int|\bar\phi_h|^2|\nabla h|^2dx=
\int_0^{H_0}|\bar\phi_h|^2
\left(\int_0^{2\pi}\frac{|\nabla \Psi|^2}{ J}d\theta\right)dh.
\end{eqnarray*}
Here $J=\Psi_y\Theta_x-\Psi_x\Theta_y$ is the Jacobian of the
coordinate change.
Note that $\nabla\Theta=\rho\nabla^\perp\Psi$ with some scalar function
$\rho>0$, so that
\[
J=\rho|\nabla \Psi|^2,~~|\nabla\Theta|=\rho|\nabla\Psi|\hbox{ and }
dl=d\theta/|\nabla\Theta|.
\]
Therefore, we have
\begin{equation}\label{jul-ah}
 \int_0^{2\pi}\frac{|\nabla\Psi|^2}{J}d\theta
=\oint_{\Psi(x,y)=h}|\nabla\Psi|dl=p(h),
\end{equation}
and thus we have a weighted $H^1$-bound
\[
\int_0^{H_0}p(h)|\bar\phi_h|^2dh<+\infty,
\]
which follows from (\ref{fr-uniform}), and hence $\bar\phi(h)$ is
continuous for $h<H_0$, as $p(h)\sim C(H_0-h)$ for $h$ close to $H_0$.

 Next, we re-write (\ref{fr-1-n}) in the
curvilinear coordinates:
\begin{eqnarray}\label{fr-2-n}
&&-\frac{|\nabla\Psi|^2}{J}\frac{\partial^2\phi_n}{\partial h^2}
-\frac{|\nabla\Theta|^2}{J}\frac{\partial^2\phi_n}{\partial\theta^2}-
\frac{(\Delta\Psi)}{J}\pdr{\phi_n}{h} -
\frac{(\Delta\Theta)}{J}\pdr{\phi_n}{\theta}+A_n\pdr{\phi_n}{\theta}
=\frac{1}{J}
\lambda g(\phi_n),\\
&&\phi_n(H_0,\theta)=0,
\hbox{ $\phi_n(h,\theta)$ is bounded for $0\le h\le H_0$.}\nonumber
\end{eqnarray}
Integrating this equation in $\theta$ and passing to the limit $n\to
+\infty$ we obtain the limit problem for the function $\bar\phi$:
\begin{eqnarray}\label{fr-2-n1}
&&-a(h)\bar\phi''(h)
-b(h)\bar\phi'(h)
=\lambda c(h)g(\bar\phi),\\
&&\bar\phi(H_0)=0,~~\hbox{ $\bar\phi(h)$ is bounded
for $0\le h\le H_0$,}\nonumber
\end{eqnarray}
with
\[
a(h)=\int_0^{2\pi}\frac{|\nabla\Psi|^2}{J}d\theta,~~~
b(h)=\int_0^{2\pi}\frac{\Delta\Psi}{J}d\theta,~~~
c(h)=\int_0^{2\pi}\frac{d\theta}{J}.
\]
It remains only to observe that (\ref{fr-2-n1}) is nothing but the
effective problem (\ref{cf-fr-3}). Indeed,  as in (\ref{jul-ah})
we compute that
\[
c(h)=\int_0^{2\pi}\frac{d\theta}{J}
=\oint_{\Psi(x,y)=h}\frac{dl}{|\nabla\Psi|}=T(h),
\]
and
\[
b(h)=\int_0^{2\pi}\frac{\Delta\Psi}{J}d\theta
=\oint_{\Psi(x,y)=h}\frac{\Delta\Psi}{|\nabla\Psi|}dl=p'(h).
\]
The last equality above follows from the fact that
\[
p(h)=\oint_{\Psi(x,y)=h}|\nabla\Psi|dl=\int_{G_h}\Delta\Psi dxdy.
\]
Here $G_h=\{h\le\Psi(x,y)\le H_0\}$ is the interior of the streamline
$\{\Psi(x,y)=h\}$. It follows that (\ref{fr-2-n1}) is, indeed,
the effective problem (\ref{cf-fr-3}), so that $\bar\phi(h)$ is a
positive solution of (\ref{cf-fr-3}). Therefore, in particular, we have
$\lambda\le\bar\lambda^*$ and (\ref{fr-lambda-lambda-bar}) holds.

\subsubsection*{ Subsolution: the case with no saddles on $\partial\Omega$}

We now prove that
\begin{equation}\label{fr-lambda-bar-lambda}
\liminf_{A\to\infty}\lambda^*(A)\ge \bar\lambda^*.
\end{equation}
Together with (\ref{fr-lambda-lambda-bar}) this will complete the
proof of Proposition~\ref{cf-fr-thm-onecell}. This is done as follows:
we take any $\lambda_0<\bar\lambda^*$ and show that
$\lambda_0\le\lambda^*(A)$ for a sufficiently large $A$ by
constructing a bounded positive super-solution to (\ref{fr-1}) with
$\lambda=\lambda_0$. However, the singular points of $\Psi(x,y)$ cause
technical difficulties in the construction of the sub-solution. Hence,
we first consider the special case when $\Psi(x,y)$ has no saddles on
$\partial\Omega$.

Let $\lambda_0<\bar\lambda^*$ and let $\bar\phi(h)$ be the
corresponding positive solution of (\ref{cf-fr-3}) with some
$\lambda\in(\lambda_0,\bar\lambda^*)$. An important observation is
that there exists $C<+\infty$ so that for all $(x,y)\neq(x_0,y_0)$
(that is, not the maximum of the stream-function $\Psi(x,y)$), we have
\begin{equation}\label{cf-laplace-bd}
|\Delta_{x,y}\bar\phi(h(x,y))|\le C,~~(x,y)\neq (x_0,y_0).
\end{equation}
To see that we write
\[
\Delta_{x,y}\bar\phi={|\nabla\Psi|^2} \bar\phi''(h(x,y))
+(\Delta\Psi)\bar\phi'(h(x,y))
\]
and note that for $h$ close to $H_0$ we have $|\nabla\Psi|^2\sim
(H_0-h)$ and $ \bar\phi''(h)\sim 1/(H-h_0)$ so that the first term
above is uniformly bounded in $(x,y)\in\Omega$.

We look for a super-solution to (\ref{fr-1}) in the form
\[
\phi=\bar\phi+\eta_A,
\]
where $\eta_A$ is smooth and bounded.
Then the uniform bound (\ref{cf-laplace-bd})
and  similar bounds
on other second derivatives of $\bar\phi(x,y)$ show that if  we have
\begin{equation}\label{cf-phi-eta-pos}
-\Delta\phi+Au\cdot\nabla\phi\ge\lambda_0 g(\phi),~~(x,y)\neq (x_0,y_0),
\end{equation}
in $\Omega$ and $\phi\ge 0$ on $\partial\Omega$, then $\phi(x,y)$ is a
weak super-solution to (\ref{fr-1}) with $\lambda=\lambda_0$ in the sense of~\cite{BCMR} and
thus $\lambda_0\le\lambda^*(A)$.
We choose $\eta_A=\eta_A(x,y)$ as the solution of
\begin{eqnarray}\label{cf-eta-eq}
&&-\Delta\eta+Au\cdot\nabla\eta=
\Delta_{x,y}\bar\phi+\lambda g(\bar\phi),\\
&&\eta(x)=0~~\hbox{ on $\partial \Omega$.}\nonumber
\end{eqnarray}
Assume that we can show that
\begin{equation}\label{cf-eta-unif-bd}
||\eta_A ||_{L^{\infty}} \to 0, \hbox{ as } A \to \infty,
\end{equation}
then~(\ref{cf-phi-eta-pos}) holds:
\begin{eqnarray*}
-\Delta\phi+Au\cdot\nabla\phi=
-\Delta\bar\phi-\Delta\eta_A+A u\cdot\nabla\eta_A=\lambda g(\bar\phi)
\ge \lambda_0 g\left(\bar\phi+\eta_A\right)=\lambda_0g(\phi),
\end{eqnarray*}
for $A$ sufficiently large.
Above we used the fact that
$\lambda_0<\lambda$ and a uniform bound for $\bar\phi$. Hence,
$\lambda_0\le\lambda^*(A)$ for $A$ sufficiently large. This will prove
Proposition~\ref{cf-fr-thm-onecell} in the special case when the
domain $\Omega$ consists of one cell and the boundary $\partial\Omega$
contains no saddles of the stream-function $\Psi$.

It remains to establish (\ref{cf-eta-unif-bd}). To this end consider a
cut-off function $\chi(s)$ so that $0\le \chi(s)\le 1$ and $\chi(s)=1$
for $|s|<1/2$ and $\chi(s)=0$ for $|s|>1$ and split
\[
\Delta_{x,y}\bar\phi+\lambda g(\bar\phi)=q_1+q_2,
\]
with
\[
q_1=\left(\Delta_{x,y}\bar\phi+\lambda g(\bar\phi)\right)
\chi\left(\frac{H_0-h(x,y))}{\delta}\right),~~
q_2=\left(\Delta_{x,y}\bar\phi+\lambda g(\bar\phi)\right)
\left[1-\chi\left(\frac{H_0-h(x,y))}{\delta}\right)\right].
\]
The small parameter $\delta>0$ is to be chosen below.
We define, accordingly, the functions $\eta_j$, $j=1,2$ as solutions of
\begin{eqnarray}\label{cf-eta12-eq}
&&-\Delta\eta_j+A u\cdot\nabla\eta_j=q_j,\\
&&\eta_j(x)=0~~\hbox{ on $\partial \Omega$,}\nonumber
\end{eqnarray}
so that $\eta_A=\eta_1+\eta_2$. As $q_1$ is uniformly bounded in
$L^\infty(\Omega)$, the function $\eta_1$ can be bounded using
Lemma~\ref{lem-lp-linfty} as
\[
\|\eta_1\|_\infty\le {C}\delta^\alpha
\]
with some $\alpha\in(0,1)$. We may further split the function
$\eta_2=\eta_3+\eta_4$, with the function $\eta_3$ that solves the ODE
along the closed streamlines:
\[
A u\cdot\nabla\eta_3=q_2,~~\eta_3(\theta=0,h)=0.
\]
This equation is solvable because
\[
\oint_{L}\left[\Delta_{x,y}\bar h+\lambda g(\bar h)\right]dl=0,
\]
as this is how the Freidlin problem is obtained. The function $\eta_3$
satisfies the estimate
\begin{equation}\label{abcdefg}
\|\eta_3\|_{C^2(\Omega)} \le\frac{F_1(\delta)}{A},
\end{equation}
with some function $F_1(\delta)$ (which may tend to infinity as $\delta\downarrow 0$).
 Finally, $\eta_4$ satisfies
\begin{eqnarray}\label{cf-eta4-eq}
&&-\Delta\eta_4+ A u\cdot\nabla\eta_4=\Delta\eta_3,\\
&&\eta_4(x)=-\eta_3(x)~~\hbox{ on $\partial \Omega$.}\nonumber
\end{eqnarray}
Once again, Lemma~\ref{lem-lp-linfty} together with the $C^2$ estimate~(\ref{abcdefg}) on
$\eta_3$ above implies that
\[
\|\eta_4\|_\infty\le\frac{CF_1(\delta)}{A}.
\]
Altogether we see that for any $\eps>0$ we can find find $\delta>0$, and then find $A_0$ so that
for any $A>A_0$
\begin{equation}\label{cf-eta-bds}
\|\eta_A\|_\infty\le \eps.
\end{equation}
This proves Proposition~\ref{cf-fr-thm-onecell} in the special case
when the domain $\Omega$ consists of one cell and the boundary
$\partial\Omega$ contains no saddles of the stream-function $\Psi$.

\subsubsection*{Approximation on a smaller domain}

Now, we establish the claim of Proposition~\ref{cf-fr-thm-onecell} for
domains $\Omega$ which consist of one cell but may have saddles of the
stream-function on the boundary $\partial\Omega$.  In order to avoid
dealing with the singular points on the boundary of $\Omega$ in the
construction a sub-solution we need to consider a slightly smaller
domain $\Omega_\eps=\{\eps\le\Psi(x,y)\le H_0\}\subset\Omega$, with
$\eps>0$ small. The domain $\Omega_\eps$ has no saddles on
$\partial\Omega$ and thus the conclusion of
Proposition~\ref{cf-fr-thm-onecell} holds for $\Omega_\eps$ by what we
have shown above.

Define $\lambda_\eps^*(A)$ as the explosion threshold for the
problem in $\Omega_\eps$:
\begin{eqnarray}\label{fr-1-eps}
&&-\Delta\phi_A^\eps+Au\cdot\nabla\phi_A^\eps=\lambda f(\phi_A^\eps)
\hbox{ in $\Omega_\eps$} \\
&&\phi_A^\eps=0\hbox{ on $\partial\Omega_\eps$.} \nonumber
\end{eqnarray}
We also let $\bar\lambda_\eps^*$ be the explosion threshold for the
corresponding Freidlin problem:
\begin{eqnarray}\label{cf-fr-3-eps}
&&-\frac{1}{T(h)}\frac{d}{dh}\left(p(h)\frac{d\bar\phi_\eps}{dh}\right)=
\lambda g(\bar\phi_\eps)\\
&&\hbox{ $\bar\phi_\eps(\eps)=0$, $\bar\phi_\eps'(h)$ is bounded for
$\eps\le h\le H_0$,}
\nonumber
\end{eqnarray}
with $T(h)$ and $p(h)$ still given by (\ref{cf-fr-one-coeff}). As we
have mentioned, since $\Omega_\eps$ has no saddles of $\Psi$ on
$\partial\Omega_\eps$, we have
\begin{equation}\label{cf-leps-barleps}
\lim_{A\to+\infty}\lambda_\eps^*(A)=\bar\lambda_\eps,
\end{equation}
for all $\eps>0$,
according to the what we have already shown above.

It is clear from the definition of $\lambda_\eps^*(A)$ and
$\bar\lambda_\eps^*$ that $\lambda_\eps^*(A)\ge\lambda^*(A)$ and
$\bar\lambda_\eps^*\ge\bar\lambda^*$. The next two lemmas show that
the passage to the limit $\eps\to 0$ is harmless.  The first statement
concerns the Freidlin thresholds.
\begin{lemma}\label{lem-barl-eps}
We have $\disp\lim_{\eps\to 0}\bar\lambda_\eps^*=\bar\lambda^*$, where
$\bar\lambda^*$ and $\bar\lambda_\eps^*$ are the explosion thresholds
of (\ref{cf-fr-3}) and (\ref{cf-fr-3-eps}).
\end{lemma}
{\bf Proof.} The proof of this lemma is rather straightforward.  It is
clear that $\bar\lambda_\eps^*\ge\bar\lambda^*$ for all $\eps>0$. On
the other hand, given $\gamma\in(0,1)$, for
$\lambda<(1-\gamma)\bar\lambda_\eps^*$ we can find $\delta>0$ which
depends only on $\gamma$ so that solution of the following problem
exists:
\begin{eqnarray}\label{cf-fr-3-eps-delta}
&&-\frac{1}{T(h)}\frac{d}{dh}\left(p(h)\frac{d\bar\phi}{dh}\right)=
\lambda g(\bar\phi)\\
&&\hbox{ $\bar\phi(\eps)=\delta$, $\bar\phi'(h)$ is bounded for
$\eps\le h\le H_0$,}
\nonumber
\end{eqnarray}
for all $\eps>0$. Then it is easy to verify that, for a sufficiently small $\eps>0$ (and
$\delta>0$ fixed),
solutions of the iteration
process
\begin{eqnarray}\label{cf-fr-3-eps-delta-bar}
&&-\frac{1}{T(h)}\frac{d}{dh}\left(p(h)\frac{d\bar\phi_n}{dh}\right)=
\lambda g(\bar\phi_{n-1})\\
&&\hbox{ $\bar\phi_n(0)=0$, $\bar\phi_n'(h)$ is bounded for
$0\le h\le H_0$,}
\nonumber
\end{eqnarray}
with $\phi_0=0$ are increasing in $n$, and satisfy
$\bar\phi_n(h)\le\delta$ for $0\le h\le\eps$ and
$\bar\phi_n(x)\le\bar\phi$ for $\eps\le h\le H_0$.  Thus, the sequence
$\bar\phi_n(h)$ converges as $n\to+\infty$ to a bounded solution of
(\ref{cf-fr-3}) so that $\lambda\le\bar\lambda$. $\Box$

The next lemma shows that $\lambda_\eps^*(A)$ is close to
$\lambda^*(A)$ for $\eps$ small.
\begin{lemma}\label{lem-leps-l0}
For any $\gamma>0$ there exists $\eps_0>0$ and $A_0$ so that for all
$\eps\in(0,\eps_0)$ and $A>A_0$ we have
$(1-\gamma/4)\lambda_\eps^*(A)\le \lambda^*(A)$.
\end{lemma}
This lemma is sufficient to show that
\begin{equation}\label{cf-lambdabar-l}
\bar\lambda^*\le\liminf_{A\to+\infty}\lambda^*(A)
\end{equation}
and thus finish the proof of Proposition~\ref{cf-fr-thm-onecell} for
all one-cell domains, as we have already established
(\ref{fr-lambda-lambda-bar}).  In order to see that
(\ref{cf-lambdabar-l}) holds, take $\gamma\in(0,1)$ and find $\eps_0$
and $A_0$ as in Lemma~\ref{lem-leps-l0}. Consider any
$\lambda<(1-\gamma)\bar\lambda^*$.  Then Lemma~\ref{lem-barl-eps}
implies that there exists $\eps_1<\eps_0$ so that
$\lambda<(1-\gamma/2)\bar\lambda_{\eps_1}^*$.  Now,
(\ref{cf-leps-barleps}) implies that we can find $A_1$ so that
$\lambda<(1-\gamma/4)\lambda_{\eps_1}^*(A)$ for all $A>A_1$. As
$\eps_1<\eps_0$ we may use Lemma~\ref{lem-leps-l0} to conclude that
$\lambda<\lambda^*(A)$ for all $A>A_0+A_1$. Therefore,
(\ref{cf-lambdabar-l}) holds. This finishes the proof of
Proposition~\ref{cf-fr-thm-onecell}. $\Box$

\subsubsection*{The proof of Lemma~\ref{lem-leps-l0}}

The proof of Lemma~\ref{lem-leps-l0} is based on the iteration
argument and stratification Lemma~\ref{cf-lem-exit}. We start with
$\lambda<(1-\gamma)\lambda_\eps^*(A)$ for all $A\geq A_0$, $\eps \leq \eps_0$ and construct a solution of the
explosion problem on the whole domain $\Omega$ by the iteration procedure.
Set $\phi_0=0$ and let $\phi_n$ be the solution of
\begin{eqnarray}\label{cf-out-7}
&&-\Delta\phi_n+Au\cdot\nabla\phi_n=\lambda g(\phi_{n-1})
\hbox{ in $\Omega$,}\\
&&\phi_n=0\hbox{ on $\partial\Omega$.}\nonumber
\end{eqnarray}
We claim that the sequence $\phi_n(x)$ is increasing in $n$, pointwise
in $x$, for each $A$ and there exists $A_0$ so that
for $A>A_0$ it has a uniformly bounded limit
$\bar\phi(x)\in L^\infty({\Omega})$ which satisfies
(\ref{fr-1}).

Pointwise monotonicity of $\phi_n(x)$ in $n$ is standard:
the difference $\eta_1(x)=\phi_1(x)-\phi_0(x)=\phi_1(x)$ satisfies
\begin{eqnarray}\label{cf-out-9}
&&-\Delta\eta_1+Au\cdot\nabla\eta_1=\lambda g(0)>0
\hbox{ in $\Omega$,}\\
&&\eta_1=0\hbox{ on $\partial\Omega$.}\nonumber
\end{eqnarray}
Hence, we have $\eta_1(x)>0$ in $\Omega$ and $\phi_1(x)>\phi_0(x)$.
Let us assume that
$\eta_n(x)=\phi_n(x)-\phi_{n-1}(x)\ge 0$ in $\Omega$. As the nonlinearity
$g(s)$ is increasing in $s$, the function
$\eta_{n+1}(x)$ is the solution of
\begin{eqnarray}\label{cf-out-10}
&&-\Delta\eta_{n+1}+Au\cdot\nabla\eta_{n+1}=
\lambda [g(\phi_n)-g(\phi_{n-1})]\ge 0
\hbox{ in $\Omega$,}\\
&&\eta_1=0\hbox{ on $\partial\Omega$.}\nonumber
\end{eqnarray}
Therefore, we have $\eta_{n+1}(x)> 0$ inside $\Omega$ and thus the
sequence $\phi_n(x)$ is increasing in $n$.

We need to show that the sequence $\phi_n(x)$ is uniformly bounded
from above.  We recall that $\lambda<(1-\gamma)\lambda_\eps^*(A)$ and
thus the minimal positive solution of (\ref{fr-1-eps}) satisfies
$0\le\phi_A^\eps\le C(\gamma)$ with the constant $C(\gamma)$
independent $\eps$ and $A$.  Hence, we may find $\delta>0$ which
depends only on $\gamma>0$ but not on $\eps$ or $A$ so that solution
of the following problem exists:
\begin{eqnarray}\label{cf-out-8}
&&-\Delta \zeta+Au\cdot\nabla \zeta=\lambda g(\zeta)
\hbox{ in $\Omega_\eps$,}\\
&&\zeta=\delta \hbox{ on $\partial \Omega_\eps.$}\nonumber
\end{eqnarray}
To see that such $\delta>0$ exists, let $\psi(x)$ satisfy
(\ref{fr-1-eps}) with $\lambda'=(1-\gamma/2)\lambda_\eps^*(A)$ and
set
$r(x)=\psi(x)+\delta$, then $r(x)$ satisfies
\begin{eqnarray}\label{cf-out-ind-delta-r}
&&-\Delta r+Au\cdot\nabla r=
\lambda' g(r-\delta)
\hbox{ in $\Omega_\eps$}\\
&&r=\delta\hbox{ on $\partial\Omega_\eps$.}\nonumber
\end{eqnarray}
It is a super-solution for (\ref{cf-out-8}) if we ensure that
$\lambda' g(r-\delta)>\lambda g(r)$, or, equivalently,
$\delta\in(0,1)$ is taken so small that we have for all $x\in\Omega$:
\begin{equation}\label{cf-g-g'}
\frac{g(\psi(x))}{g(\psi(x)+\delta)}\ge\frac{1-\gamma}{1-\gamma/2}.
\end{equation}
The function $\psi(x)$ obeys a uniform bound
$\|\psi\|_{L^\infty(\Omega_\eps)}\le K(\gamma)$ with $K(\gamma)$
independent of $\eps$ and $A$.  Let $M=\sup_{0\le s\le K(\gamma)+1}
g'(s)$, then (\ref{cf-g-g'}) is guaranteed if we have
\[
\frac{g(s)}{g(s)+M\delta}\ge\frac{1-\gamma}{1-\gamma/2}
\]
for all $s\in[0,K(\gamma)]$. A direct computation shows that this is
possible if we take
\[
\delta<\delta_0=\frac{\gamma g(0)}{M(1-\gamma)}.
\]
Under this assumption $r(x)$ provides a super-solution for
(\ref{cf-out-8}) and thus a positive solution of this problem can
be constructed by the standard iteration procedure.

In order to show that the sequence $\phi_m(x)$ is bounded we use $\zeta(x)$, the minimal positive solution of
(\ref{cf-out-8}) and we need the
following analog of the stratification Lemma~\ref{cf-lem-exit}.
\begin{lemma}\label{cf-lem-omepsbd}
Fix $M>0$. Let $q_A(x)$ be solution of
\begin{eqnarray}\label{cf-eq-omeps}
&&-\Delta q_A+Au\cdot\nabla q_A=M\hbox{ in $\Omega$},\\
 &&q_A=0\hbox{ on $\partial\Omega$.}\nonumber
\end{eqnarray}
Then, given any $\delta>0$ there exist $A_0>0$ and $\eps_0>0$ so that
for all solutions of (\ref{cf-eq-omeps}) we have $0\le
q_A(x)\le\delta$ in $G_\eps=\Omega\setminus\Omega_\eps$ for all $A>A_0$
and $\eps<\eps_0$.
\end{lemma}
We do not present the proof of this lemma as it is essentially
contained in that of Lemma~\ref{cf-lem-exit}.

We choose $M>0$ so that $\lambda g(\zeta(x))\le M$
for all  $x\in\Omega_\eps$, where  $\zeta(x)$ is the minimal positive solution of (\ref{cf-out-8}).
We may also take $A>0$ sufficiently large, as in Lemma~\ref{cf-lem-omepsbd}. We claim that then we will have, for all $m\ge
1$, (i) $0\le\phi_m(x)\le\delta$ in $G_\eps=\Omega\setminus\Omega_\eps$,
and (ii) $0\le\phi_m(x)\le \zeta(x)$ for all $x\in\Omega_\eps$.

Let us prove this by induction. The function $\phi_1(x)$ satisfies
\begin{eqnarray}\label{cf-out-11}
&&-\Delta\phi_1+Au\cdot\nabla\phi_1=\lambda g(0)
\hbox{ in $\Omega$}\\
&&\phi_1=0\hbox{ on $\partial\Omega$.}\nonumber
\end{eqnarray}
Our choice of $M$ ensures that the right side in (\ref{cf-out-11}) is
bounded above by $M$. Thus, if the amplitude $A$ is sufficiently large
we have $0\le\phi_1(x)\le\delta$ in the tube $G_\eps$ and in
particular on $\partial\Omega_\eps$. Therefore, the difference
$s_1(x)=\zeta(x)-\phi_1(x)$ satisfies
\begin{eqnarray}\label{cf-out-11a}
&&-\Delta s_1+Au\cdot\nabla s_1=\lambda [g(q_j)-g(0)]\ge 0
\hbox{ in $\Omega_\eps$},\\
&&s_1\ge 0\hbox{ on $\partial\Omega_\eps$.}\nonumber
\end{eqnarray}
It follows that $s_1(x)\ge 0$ and $0\le \phi_1(x)\le\zeta(x)$ in
$\Omega_\eps$ so that our claim holds for $n=1$. Assume now
that (i) and (ii) hold for $\phi_{n-1}(x)$. It follows from the
induction assumption that
\[
\lambda g(\phi_{n-1}(x))\le \lambda g(\zeta(x))\le
M\hbox{ in $\Omega_\eps$},
\]
and $\lambda g(\phi_{n-1}(x))\le \lambda g(\delta)\le M$ in $G_\eps$.
Therefore, for $\phi_n(x)$ we have
\begin{eqnarray}\label{cf-out-12}
&&-\Delta\phi_n+Au\cdot\nabla\phi_n
=\lambda g(\phi_{n-1})\le M,
\hbox{ in $\Omega$}\\
&&\phi_n=0\hbox{ on $\partial\Omega$.}\nonumber
\end{eqnarray}
Hence, by Lemma~\ref{cf-lem-omepsbd} we have $0\le\phi_m(x)\le\delta$
in $G_\eps$. On the other hand, the difference
$s_m(x)=\zeta(x)-\phi_m(x)$ inside ${\partial\Omega_\eps}$ obeys
\begin{eqnarray}\label{cf-out-13}
&&-\Delta s_m+Au\cdot\nabla s_m=
\lambda [g(\zeta)-g(\phi_{m-1})]\ge 0
\hbox{ in $\Omega_\eps$},\\
&&s_m\ge 0\hbox{ on $\partial\Omega_\eps$,}\nonumber
\end{eqnarray}
and thus $s_m(x)\ge 0$, so that $0\le \phi_m(x)\le \zeta(x)$ in
$\Omega_\eps$.  Therefore, the sequence $\phi_m$ is increasing and
uniformly bounded from above. The limit $\bar\phi$ is a positive
solution of (\ref{fr-1}) and thus $\lambda\le\lambda^*(A)$.
This finishes the proof of Lemma~\ref{lem-leps-l0}. $\Box$

\subsection{General two-dimensional cellular flows}\label{sec:cell}

We now look at the explosion problem
\begin{eqnarray}\label{out-1}
&&-\Delta\phi^A+Au\cdot\nabla\phi^A=\lambda g(\phi^A)
\hbox{ in $\Omega$,}\\
&&\phi^A=0\hbox{ on $\partial\Omega$,}\nonumber
\end{eqnarray}
in a two-dimensional domain $\Omega\subset\Rm^2$ with a cellular flow
$u$ which may contain more than one cell in $\Omega$ and complete the
proof of Theorem~\ref{cf-fr-thm}.
Let $\lambda_j^*(A)$ be the explosion
threshold for the problem inside each cell:
\begin{eqnarray}\label{out-1j}
&&-\Delta\phi+Au\cdot\nabla\phi=\lambda g(\phi)
\hbox{ in ${\cal C}_j$,}\\
&&\phi^A=0\hbox{ on ${\cal C}_j$.}\nonumber
\end{eqnarray}
We already know that
\[
\lim_{A\to+\infty}\lambda_j^*(A)=\bar\lambda_j^*,
\]
from Proposition~\ref{cf-fr-thm-onecell} and, of course,
$\lambda_j^*(A)\ge\lambda^*(A)$ for all $j$. Hence, all we need to
verify for the proof of Theorem~\ref{cf-fr-thm} is that for any
$\lambda<\lim_{A\to+\infty}\lambda_j^*(A)$ solution of the problem
(\ref{out-1}) on the whole domain $\Omega$ exists.

The proof is quite similar to the last part of the proof of
Proposition~\ref{cf-fr-thm-onecell}: we construct the solution of
(\ref{out-1}) by the iteration procedure. Set
\[
\lambda_0=\min_j\left[\lim_{A\to+\infty}\lambda_j^*(A)\right]
\]
and take $\gamma>0$ fixed.  Consider any
$\lambda\in(0,(1-\gamma)\lambda_0)$, start the iteration process
with $\phi_0=0$ and define $\phi_m$ as the solution of
\begin{eqnarray}\label{out-7}
&&-\Delta\phi_m^A+Au\cdot\nabla\phi_m^A=\lambda g(\phi_{m-1}^A)
\hbox{ in $\Omega$,}\\
&&\phi_m^A=0\hbox{ on $\partial\Omega$.}\nonumber
\end{eqnarray}
\begin{lemma}\label{lem-iteration}
The sequence $\phi_m^A(x)$ is increasing in $m$, pointwise in $x$, for
each $A>A_0$ and has a uniformly bounded limit $\bar\phi_A(x)\in
L^\infty({\Omega})$ which is a solution of (\ref{out-1}).
\end{lemma}

\subsubsection*{The proof of Lemma~\ref{lem-iteration}.}

The proof is quite analogous to that of Lemma~\ref{lem-leps-l0}
except that we use Lemma~\ref{cf-lem-exit} where Lemma~\ref{cf-lem-omepsbd}
was used before.  The
first increment $\eta_1(x)=\phi_1^A(x)-\phi_0^A(x)=\phi_1^A(x)$
satisfies
\begin{eqnarray}\label{out-9}
&&-\Delta\eta_1+Au\cdot\nabla\eta_1=\lambda g(0)>0
\hbox{ in $\Omega$,}\\
&&\eta_1=0\hbox{ on $\partial\Omega$.}\nonumber
\end{eqnarray}
Hence, we have $\eta_1(x)>0$ in $\Omega$ and
$\phi_1(x)>\phi_0(x)$. Let us assume that
$\eta_n(x)=\phi_n^A(x)-\phi_{n-1}^A(x)\ge 0$ in $\Omega$. As the
nonlinearity $g(s)$ is increasing in $s$, the function $\eta_{n+1}(x)$
is the solution of
\begin{eqnarray}\label{out-10}
&&-\Delta\eta_{n+1}+Au\cdot\nabla\eta_{n+1}=
\lambda [g(\phi_n)-g(\phi_{n-1})]\ge 0
\hbox{ in $\Omega$},\\
&&\eta_1=0\hbox{ on $\partial\Omega$.}\nonumber
\end{eqnarray}
Therefore, we have $\eta_{n+1}(x)> 0$ inside $\Omega$ and thus the
sequence $\phi_n^A(x)$ is increasing in $n$. We need to show that it
is uniformly bounded from above.
As in the proof of Lemma~\ref{lem-leps-l0}, we can find $\delta>0$
so that solution of
the following explosion problem exists on each cell ${\cal C}_j$ for
all $A>A_0$, $\delta\in(0,\delta_0)$ and $\lambda<(1-\gamma)\lambda_0$:
\begin{eqnarray}\label{out-ind-delta}
&&-\Delta q_j+Au\cdot\nabla q_j=
\lambda g(q_j)
\hbox{ in ${\cal C}_j$}\\
&&q_j=\delta\hbox{ on $\partial{\cal C}_j$,}\nonumber
\end{eqnarray}
and, moreover, $0\le q_j(x)\le K(\gamma)$ with the
constant $K(\gamma)$ which does not depend on the flow amplitude $A$.
This is shown by exactly the same argument we used to show the existence
of a positive solution for (\ref{cf-out-8}).


As $q_j(x)$ are uniformly bounded by $K_j(\gamma)$, we may choose
$M>0$ so that $\lambda g(q_j(x))\le M$ for all $j$ and all
$x\in\Omega$.  We also take $A>0$ sufficiently large, as in
Lemma~\ref{cf-lem-exit} (but with the right side of (\ref{cf-out-14})
replaced by the constant $M$ rather than $1$). We claim that the
following bounds will be preserved by the iteration procedure: for all
$m\ge 1$, (i) $0\le\phi_m(x)\le\delta$ on the skeleton of separatrices
${\cal D}_0$,
and (ii) $0\le\phi_m(x)\le q_j(x)$ for all $x\in{\cal C}_j$.
Again, we check this by induction. The function $\phi_1^A(x)$ satisfies
\begin{eqnarray}\label{out-11}
&&-\Delta\phi_1^A+Au\cdot\nabla\phi_1^A=\lambda g(0),
\hbox{ in $\Omega$,}\\
&&\phi_1^A=0\hbox{ on $\partial\Omega$.}\nonumber
\end{eqnarray}
The right side in (\ref{out-11}) is bounded above by $M$ and hence,
according to Lemma~\ref{cf-lem-exit}, for the amplitude $A$
sufficiently large we have $0\le\phi_1^A(x)\le\delta$ on the skeleton
${\cal D}_0$. Therefore, the first increment $s_1(x)=q_j(x)-\phi_1(x)$
satisfies
\begin{eqnarray}\label{out-11a}
&&-\Delta s_1+Au\cdot\nabla s_1=\lambda [g(q_j)-g(0)]\ge 0
\hbox{ in ${\cal C}_j$},\\
&&s_1\ge 0\hbox{ on $\partial{\cal C}_j$,}\nonumber
\end{eqnarray}
and thus $s_1(x)\ge 0$ and $0\le \phi_1^A(x)\le q_j(x)$ in ${\cal
C}_j$. Let us now assume that (i) and (ii) are true for
$\phi_{m-1}(x)$ and show that they hold for $\phi_m(x)$ -- the
argument is exactly as for $m=1$. By the induction assumption we have
\[
\lambda g(\phi_{m-1}(x))\le \lambda g(q_j(x))
\le M\hbox{ in the cell ${\cal C}_j$},
\]
and thus $\phi_m(x)$ satisfies
\begin{eqnarray}\label{out-12}
&&-\Delta\phi_m^A+Au\cdot\nabla\phi_m^A
=\lambda g(\phi_{m-1}^A)\le \lambda g(q_j(x))\le M,
\hbox{ in $\Omega$,}\\
&&\phi_m^A=0\hbox{ on $\partial\Omega$.}\nonumber
\end{eqnarray}
Hence, by Lemma~\ref{cf-lem-exit} we have
$0\le\phi_m^\eps(x)\le\delta$ on ${\cal D}_0$. Now, the difference
$s_m(x)=q_j(x)-\phi_m(x)$ inside ${\cal C}_j$ satisfies
\begin{eqnarray}\label{out-13}
&&-\Delta s_m+A u\cdot\nabla s_m=
\lambda [g(q_j)-g(\phi_{m-1}^\eps)]\ge 0
\hbox{ in ${\cal C}_j$},\\
&&s_m\ge 0\hbox{ on $\partial{\cal C}_j$,}\nonumber
\end{eqnarray}
and thus $s_m(x)\ge 0$, so that $0\le \phi_m^A(x)\le q_j(x)$ in ${\cal
C}_j$.  Therefore, the sequence $\phi_m^A$ is increasing and uniformly
bounded from above. The limit $\bar\phi^A(x)$ is a positive solution
of (\ref{out-1}). This completes the proof of Lemma~\ref{lem-iteration}
and thus that of Theorem~\ref{cf-fr-thm}. $\Box$

\end{document}